\documentclass[10pt]{article}
\usepackage{amssymb}
\usepackage{amsmath}
\usepackage{bm}
\usepackage{amsthm}
\usepackage{amsbsy,tcolorbox,mdframed,pdfpages,makecell}
\usepackage{psfrag}
\usepackage{graphics}
\usepackage{graphicx}
\usepackage{pgfplots}
\usepackage{relsize}
\usepackage{color}
\usepackage{xcolor}
\usepackage{pdfpages}  
\definecolor{mycolor}{RGB}{220,220,220} 
\usepackage{url}
\usepackage{hyperref}
\usepackage[square,sort&compress,comma,numbers]{natbib}
\usepackage{graphicx}
\usepackage{caption}
\usepackage{subcaption}
\hypersetup{
	colorlinks=true,
	linkcolor=blue,
	filecolor=magenta,      
	urlcolor=cyan,
}
\theoremstyle{plain}
\newtheorem{theorem}{Theorem}[section]

\newtheorem{assumption}[theorem]{Assumption}

\newtheorem{lemma}[theorem]{Lemma}

\newcommand{\norm}[1]{\lVert #1 \rVert}
\newcommand{\hhonorm}[1]{{\lVert #1 \rVert}_{1,h}}
\newcommand{\ru}[1]{R_{\mathcal{T}}\widehat{#1}_h}
\newcommand{\hhoelem}[1]{\widehat{#1}_h}
\newcommand{\lochhoelem}[1]{\widehat{#1}_T}
\newtheorem{remark}[theorem]{Remark}

\author{Gouranga Mallik \thanks{Department of Mathematics, School of Advanced Sciences, Vellore Institute of Technology, Vellore - 632014, India (\texttt{gouranga.mallik@vit.ac.in})} \and Ramesh Chandra Sau \thanks{Department of Mathematics, Indian Institute of Technology, Bombay, India (\texttt{rcsau@math.iitb.ac.in, rcsau1994@gmail.com})}}

\title{An Error Analysis of Second Order Elliptic Optimal Control Problem via Hybrid Higher Order Methods}

\begin{document}
\allowdisplaybreaks[4]
\numberwithin{figure}{section}

%
\numberwithin{table}{section}
 \numberwithin{equation}{section}

\date{}

\maketitle
\allowdisplaybreaks
\def\R{\mathbb{R}}
\def\I{\mathbb{I}}
\def\dx{{\rm~dx}}
\def\dy{{\rm~dy}}
\def\ds{{\rm~ds}}
\def\y{\mathbf{y}}
\def\u{\mathbf{u}}
\def\V{\mathbf{V}}
\def\w{\mathbf{w}}
\def\J{\mathbf{J}}
\def\f{\mathbf{f}}
\def\g{\mathbf{g}}
\def\div{{\rm div~}}
\def\v{\mathbf{v}}
\def\z{\mathbf{z}}
\def\b{\mathbf{b}}
\def\cV{\mathcal{V}}
\def\Vsp{\hat{V}^k_{h,0}}
\def\rt{R_{\mathcal{T}}}
\def\ih{\widehat{I}_h}
\def\V{\mathbf{V}}
\def\cT{\mathcal{T}}
\def\sjump#1{[\hskip -1.5pt[#1]\hskip -1.5pt]}

\begin{abstract}
    This paper presents the design and analysis of a Hybrid High-Order (HHO) approximation for a distributed optimal control problem governed by the Poisson equation. We propose three distinct schemes to address unconstrained control problems and two schemes for constrained control problems. For the unconstrained control problem, while standard finite elements achieve a convergence rate of \( k+1 \) (with \( k \) representing the polynomial degree), our approach enhances this rate to \( k+2 \) by selecting the control from a carefully constructed reconstruction space. For the box-constrained problem, we demonstrate that using lowest-order elements (\( \mathbb{P}_0 \)) yields linear convergence, in contrast to finite element methods (FEM) that require linear elements to achieve comparable results. Furthermore, we derive a cubic convergence rate for control in the variational discretization scheme. Numerical experiments are provided to validate the theoretical findings.
\end{abstract}


\section{Introduction}
The study of optimal control problems governed by partial differential equations (PDEs) is a significant area of research in applied mathematics. The optimal control problem consists of finding a control variable that minimizes a cost functional subject to a PDE. 
Due to its importance in applications, several numerical methods have been proposed to approximate the solutions. The finite element approximation of the optimal control problem started with the work of Falk \cite{Falk:1973:Control} and Geveci \cite{Gevec:1979:Control}. 
A control can act in the interior of a domain, in this case, we call distributed, or on the boundary of a domain, we call boundary (Neumann or Dirichlet) control problem. We refer to \cite{Hinze:2011:Control,meyerrosch:2004:Optiaml,Sudipto:2015:IP,Ramesh:2018} for distributed control related problem, to \cite{CasasM:2008:Neum,Casas_Dhamo_2012,Sudipto:2015:IP,Ramesh:2018} for the Neumann boundary control problem, and to \cite{Hinze:2005:Control} for a variational discretization approach. The Dirichlet boundary control problem has been studied in \cite{CMR:2009:Dirich,CasasRaymond:2006:Dirich}.

\noindent
In this article, we consider the following optimal control problem without/with control constraints: 
\begin{equation}\label{eqn:conti_costfun}
 \min J(y, u):=\frac{1}{2}\|y-y_{d}\|_{L^2(\Omega)}^{2}+\frac{\lambda}{2}\|u\|_{L^2(\Omega)}^{2},   
\end{equation}
subject to PDE,
\begin{equation}\label{eqn:conti_state}
\begin{aligned}
-\Delta y & =f+u \quad \text{in } \Omega, \\
y & =0 \quad \quad \text { on } \partial \Omega,
\end{aligned}
\end{equation}
 where control $u$ from the space $L^2(\Omega)$ (i.e., without constraints case), or from a box-constrained set $U_{ad}$ (i.e., with constraints case). Our goal is to approximate the solution of the problem \eqref{eqn:conti_costfun}-\eqref{eqn:conti_state} using hybrid higher order (HHO) elements.

 HHO methods were introduced in \cite{ErnDiPietro2014} for linear diffusion and in \cite{Ern_DiPietro_elasticity2015} for locking-free linear elasticity. In such methods, discrete unknowns are attached to the mesh cells and to the mesh faces. The two key ingredients to devise HHO methods are a local reconstruction operator and a local stabilization operator in each mesh cell. HHO methods offer various attractive features, such as the support of polytopal meshes, optimal error estimates, local conservation properties, and computational efficiency due to compact stencils and local elimination of the cell unknowns by static condensation. As a result, these methods have been developed extensively over the past few years and now cover a broad range of applications; we refer the reader to the two recent monographs \cite{Book:HHO_Ern:2021,Book:JeromDiPietro2020} for an overview. As shown in \cite{CockburnErnDiPietro2016}, HHO methods can be embedded into the broad framework of hybridizable discontinuous Galerkin methods, and they can be bridged to nonconforming virtual element methods.

For the unconstrained control problem, standard finite element methods achieve a convergence rate of \( k+1 \). Our focus is on enhancing this rate of convergence. To address this, we propose three discretization schemes for unconstrained control problems.
 \begin{itemize}
 	\item In Scheme UC1 (UC abbreviates unconstrained),  we choose discrete control from piecewise $k$-degree polynomial space. Discrete state and adjoint state are chosen from $k$-degree HHO space $\Vsp$ (it consists of piecewise $k$-degree polynomials in the interior and on the skeleton). We obtain $k+1$ rate of convergence in the $L^2$ error of control and energy error of state and adjoint state.   

  \item  Scheme UC2 is popularly known as variational discretization; here, we do not discretize the control separately; we only discretize state and adjoint state. Thus we determine the optimal control by the adjoint state. For this scheme, we choose discrete state and adjoint state from $k$- degree HHO space $\Vsp$. This scheme shows the $k+1$ rate of convergence for control, state, and adjoint state variables.
 	
 	\item To improve the above rate of convergence, we propose the third scheme (Scheme UC 3.1 \& 3.2), and here we use control space as reconstruction of some HHO space $R_{\mathcal{T}}\widehat{V}^k_h$. This method improves the $L^2-$ rate of convergence of control to $k+2$ for $k\geq 0$. We split the third scheme into two cases. In the first scheme, UC3.1 (we call it full reconstruction), we choose the discrete state and adjoint state from $k$-degree HHO space $\Vsp$, and we prove $k+2$ convergence for polynomial degree $k=0,1$. In the second scheme UC3.2 (we call it partial reconstruction), we use a mixed order HHO space $\widehat{V}^{k+}_{h,0}$ for discrete state and adjoint state variables for $k\geq 2$. The space $\widehat{V}^{k+}_{h,0}$ consists of $k+1$ degree polynomial in the interior and $k$ degree polynomial on the skeleton. Thus we improve the convergence rates to $k+2,$ for $k\geq 2$.
 \end{itemize}
 
 Table \ref{table:Rates_UCM} shows the convergence rates and spaces for Schemes UC1, UC2, UC3.1, and UC3.2.
\begin{table}[h!]
\centering
\begin{tabular}{|p{2.0cm}|p{3.5cm}|p{3.0cm}|p{3cm}|p{3cm}|p{2.5cm}|}
\hline
Schemes & Spaces of state, \newline control and adjoint state variables  & $L^2-$ error \newline of control variable  & Energy error \newline of state variable & Energy error \newline of adjoint state variable \\ \hline\hline
Scheme UC1 & $(\hhoelem{y},u_h,\hhoelem{\phi})\in \newline \Vsp \times \mathbb{P}^k(\mathcal{T}_h) \times \Vsp,$ \newline $k\geq 0$ & $\norm{u-u_h}_{L^2(\Omega)} \newline = \mathcal{O}(h^{k+1})$ & $\hhonorm{\ih y -\hhoelem{y}} \newline = \mathcal{O}(h^{k+1})$  & $\hhonorm{\ih \phi -\hhoelem{\phi}} \newline = \mathcal{O}(h^{k+1})$ \\ \hline
Scheme UC2 \newline (Variational Discretization) & $(\hhoelem{y},u_h,\hhoelem{\phi})\in \newline \Vsp \times  L^2(\Omega) \times \Vsp,$ \newline $k\geq 0$ & $\norm{u-u_h}_{L^2(\Omega)} \newline = \mathcal{O}(h^{k+1})$ & $\hhonorm{\ih y -\hhoelem{y}} \newline = \mathcal{O}(h^{k+1})$  & $\hhonorm{\ih \phi -\hhoelem{\phi}} \newline = \mathcal{O}(h^{k+1})$ \\ \hline
Scheme UC3.1 \newline (full reconstruction) & $(\hhoelem{y},\ru{u},\hhoelem{\phi})\in \newline \Vsp \times R_{\mathcal{T}}\widehat{V}^k_h \times \Vsp,$ \newline $k=0,1$ & $\norm{u-\ru{u}}_{L^2(\Omega)} \newline = \mathcal{O}(h^{k+2})$ & $\hhonorm{\ih y -\hhoelem{y}} \newline = \mathcal{O}(h^{k+1})$  & $\hhonorm{\ih \phi -\hhoelem{\phi}} \newline = \mathcal{O}(h^{k+1})$ \\ \hline
Scheme UC3.2 \newline (partial reconstruction) & $(\hhoelem{y},\ru{u},\hhoelem{\phi})\in \newline \widehat{V}^{k+}_{h,0} \times R_{\mathcal{T}}\widehat{V}^{k}_h \times \widehat{V}^{k+}_{h,0},$ \newline $k\geq 2$ & $\norm{u-\ru{u}}_{L^2(\Omega)} \newline = \mathcal{O}(h^{k+2})$ & $\hhonorm{\ih y -\hhoelem{y}} \newline = \mathcal{O}(h^{k+1})$  & $\hhonorm{\ih \phi -\hhoelem{\phi}} \newline = \mathcal{O}(h^{k+1})$ \\ \hline
\end{tabular}
\caption{Convergence rates table for unconstraints control problem \eqref{eqn:conti_costfun}-\eqref{eqn:conti_state}.}
\label{table:Rates_UCM}
\end{table}

\noindent
For the with-constraints case, we describe three HHO schemes to analyze the control problem \eqref{eqn:conti_costfun}-\eqref{eqn:conti_state}. The regularity of constrained control problem is limited and in this case we have solutions $y\in H^3(\Omega), \phi \in H^3(\Omega),\; \text{and}\; u\in W^{1,p}(\Omega)$, with $2<p<\infty,$ when data $y_d,f\in H^1(\Omega).$ Therefore, we can not expect higher-order convergences. In the following, we briefly describe two proposed schemes and their convergence behavior.
\begin{itemize}
\item In Scheme WC1 (WC abbreviates with unconstrints) we choose discrete control from piecewise constant polynomial space. The discrete state and adjoint state variables are chosen from $0$-th degree HHO space $\widehat{V}^0_{h,0}$ (it consists of piecewise constant polynomials in the interior and on the boundary). We observe a linear rate of convergence for control, state, and adjoint state variables.

 \item  Scheme WC2 considers the variational discretization. As usual, in variational discretization, we do not discretize the control variables; we determine control from the discrete adjoint state. We choose discrete state and adjoint state from mixed degree HHO space $\widehat{V}^{1+}_{h,0}$. We improve the rate of convergence to the third order for the control variable and the second order for the state and adjoint state variable under the previously mentioned regularity.
\end{itemize}
Table \ref{table:WCM} shows the rates of convergence for Schemes WC1 and WC2.  
\begin{table}[h!]
\centering
\begin{tabular}{|p{2.0cm}|p{3.5cm}|p{2.5cm}|p{3cm}|p{3cm}|p{3cm}|}
\hline
Schemes & Spaces of state, \newline control and adjoint state variables  & $L^2-$ error \newline of control variable  & Energy error \newline of state variable & Energy error \newline of adjoint state variable \\ \hline\hline
Scheme WC1& $(\hhoelem{y},u_h,\hhoelem{\phi})\in \newline \widehat{V}^0_{h,0} \times \mathbb{P}^0(\mathcal{T}_h) \times \widehat{V}^0_{h,0},$ & $\norm{u-u_h}_{L^2(\Omega)} \newline  = \mathcal{O}(h)$ & $\norm{\ih y -\hhoelem{y}}_{1,h} \newline = \mathcal{O}(h)$  & $\norm{\ih \phi -\hhoelem{\phi}}_{1,h} \newline =  \mathcal{O}(h)$ \\ \hline
Scheme WC2 \newline (Variational discretization) & $(\hhoelem{y},u_h,\hhoelem{\phi})\in \newline \widehat{V}^{1+}_{h,0} \times U_{ad} \times \widehat{V}^{1+}_{h,0},$ & $\norm{u-u_h}_{L^2(\Omega)} \newline = \mathcal{O}(h^3)$ & $\norm{\ih y -\hhoelem{y}}_{1,h} \newline = \mathcal{O}(h^2)$  & $\norm{\ih \phi -\hhoelem{\phi}}_{1,h} \newline = \mathcal{O}(h^2)$ \\ \hline
\end{tabular}
\caption{Convergence rates table for control problem with control constrains \eqref{eqn:conti_costfun}-\eqref{eqn:conti_state}.}
\label{table:WCM}
\end{table}

\noindent
The following are some novelty and achievements of this paper:
\begin{itemize}
    \item \textit{Without constraints case:}
 Scheme UC1 uses $k$-degree piece-wise polynomials for control, state, and adjoint state variables. We observe that the energy estimate of state and adjoint state shows better convergence than finite element-based methods. Schemes UC3.1 \& UC3.2 show $k+2$ convergence rates for the $L^2$- reconstruction error of control. This is a phenomenal improvement for the rate of convergence of control error compared to the finite elements.
 \item  \textit{With constraints case:}
The standard lowest order finite elements for constrained optimal control problem needs $\mathbb{P}_0$ elements for the control variable and $\mathbb{P}_1$ elements for the state and adjoint state to produce linear convergence for the errors. However, in Scheme WC1, we observe that the lowest-order HHO elements, i.e., $\mathbb{P}_0$ elements ensure the linear rate of convergence. Therefore, there is a significant improvement in both the theory and computation. In Scheme WC2, we use a mixed degree HHO space $\widehat{V}^{1+}_{h,0}$ to obtain third-order convergence in control error and second-order convergence in state and adjoint state. Therefore it is showing better convergence than finite elements. 
\end{itemize}

The following is a breakdown of the rest of the article. Some preliminaries on HHO method have been introduced in Section \ref{sec:prelimin_hho}. Section \ref{sec:OCP_HHO_UC}, introduces three HHO-based schemes for the unconstrained control problem. Another three HHO-based schemes for the constrained control problem have been discussed in Section \ref{sec:OCP_HHO_WC}. Section \ref{sec:Numerical_experiments} has been devoted to numerical experiments.

\section{Preliminaries on HHO method}\label{sec:prelimin_hho}

\subsection{Discrete setting for HHO method.} Let $\Omega\subset \mathbb{R}^d$ be a convex polygonal/polyhedral domain with $d=2,3$. We consider a sequence of meshes $(\mathcal{T}_{h})_{h>0}$, where the parameter $h$ denotes the mesh size and goes to zero during the refinement process. For all $h>0$, we assume that the mesh $\mathcal{T}_{h}$ covers $\Omega$ exactly and consists of a finite collection of nonempty disjoint open polygonal/polyhedral cells $T$ such that $\bar{\Omega}=\cup_{T \in \mathcal{T}_{h}} \bar{T}$ and $h=\max _{h \in \mathcal{T}_{h}} h_{T}$, where $h_{T}$ is the diameter of $T$. A closed subset $F$ of $\Omega$ is defined to be a mesh face if it is a subset of a straight line/an affine hyperplane $H_{F}$ with positive $d-1$ dimensional measure and if either of the following two statements holds true: (i) there exist $T_{1}$ and $T_{2}$ in $\mathcal{T}_{h}$ such that $F \subset \partial T_{1} \cap \partial T_{2} \cap H_{F}$; in this case, the face $F$ is called an internal face; (ii) there exists $T \in \mathcal{T}_{h}$ such that $F \subset \partial T \cap \partial \Omega \cap H_{F}$; in this case, the face $F$ is called a boundary face. The set of mesh faces is a partition of the mesh skeleton, i.e., $\cup_{T \in \mathcal{T}_{h}} \partial T=\cup_{F \in \mathcal{F}_{h}} \bar{F}$, where $\mathcal{F}_{h}:=\mathcal{F}_{h}^{i} \cup \mathcal{F}_{h}^{b}$ is the collection of all faces that is the union of the set of all the internal faces $\mathcal{F}_{h}^{i}$ and the set of all the boundary faces $\mathcal{F}_{h}^{b}$. Let $h_{F}$ denote the diameter of $F \in \mathcal{F}_{h}$. For each $T \in \mathcal{T}_{h}$, the set $F_{T}:=\{F \in \mathcal{F}_{h} \mid F \subset \partial T\}$ denotes the collection of all faces contained in $\partial T, \boldsymbol{n}_{T}$ denotes the unit outward normal to $T$, and we set $\boldsymbol{n}_{T F}:=\boldsymbol{n}_{T}|_{F}$ for all $F \in \mathcal{F}_{h}$. Following \cite[Definition 1]{Ern_DiPietro_elasticity2015}, we assume that the mesh sequence $(\mathcal{T}_{h})_{h>0}$ is admissible in the sense that, for all $h>0, \mathcal{T}_{h}$ admits a matching simplicial sub-mesh $\mathcal{T}_{h}$ (i.e., every cell and face of $\mathcal{T}_{h}$ is a subset of a cell and a face of $\mathcal{T}_{h}$, respectively) so that the mesh sequence $(\mathcal{T}_{h})_{h>0}$ is shape-regular in the usual sense and all the cells and faces of $\mathcal{T}_{h}$ have a uniformly comparable diameter to the cell and face of $\mathcal{T}_{h}$ to which they belong. Owing to \cite[Lemma 1.42]{Ern_DiPietro_Book2012}, for $T \in \mathcal{T}_{h}$ and $F \in \mathcal{F}_{T}, h_{F}$ is comparable to $h_{T}$ in the sense that

$$
\varrho^{2} h_{T} \leq h_{F} \leq h_{T},
$$
where $\varrho$ is the mesh regularity parameter. Moreover, there exists an integer $N_{\partial}$ depending on $\varrho$ and $d$ such that (see \cite[Lemma 1.41]{Ern_DiPietro_Book2012} ):
$$
\max _{T \in \mathcal{T}_{h}} \operatorname{card}(\mathcal{F}_{T}) \leq N_{\partial}.
$$

\subsection{Discrete spaces} Let $k \geq 0$ be the polynomial degree. For  all $T \in \mathcal{T}_{h}$, the local space of degrees of freedom (DOFs) is defined as:
\begin{equation*}
\widehat{V}_{T}^{k}:=\mathbb{P}_{d}^{k}(T) \times\big\{\underset{F \in \mathcal{F}_{T}}{\times} \mathbb{P}_{d-1}^{k}(F)\big\}, 
\end{equation*}
where $\mathbb{P}_{d}^{k}(T)$ (resp. $\mathbb{P}_{d-1}^{k}(F)$) is the space of polynomials of degree at most $k$ on $T \in \mathcal{T}_{h}$ (resp. $F \in \mathcal{F}_{h}$). The space of piecewise polynomial of degree $k$ on $\mathcal{T}_h$ is denoted by $\mathbb{P}^k_{d}(\mathcal{T}_h)$.

\medskip \noindent
Patch all the local spaces to define the global space of DOFs as:
$$
\widehat{V}_{h}^{k}:=\big\{\underset{T \in \mathcal{T}_{h}}{\times} \mathbb{P}_{d}^{k}(T)\big\} \times\big\{\underset{F \in \mathcal{F}_{h}}{\times} \mathbb{P}_{d-1}^{k}(F)\big\}.
$$
Impose the zero boundary condition on the above global space to introduce
$$
\widehat{V}_{h, 0}^{k}:=\{\widehat{v}_{h}=((v_{T})_{T \in \mathcal{T}_{h}},(v_{F})_{F \in \mathcal{F}_{h}}) \in \widehat{V}_{h}^{k} \mid v_{F} \equiv 0 \quad \forall F \in \mathcal{F}_{h}^{b}\}.
$$
For $\widehat{v}_{h}:=(v_{\mathcal{T}},v_{\partial\mathcal{T}})\in \widehat{V}_{h, 0}^{k}$, where $v_{\mathcal{T}}\in \mathbb{P}^k_d(\mathcal{T}_h)$ as $v_{\mathcal{T}}|_T:=v_T$ and $v_{\partial\mathcal{T}}|_F:=v_F$.
Moreover, a mixed order HHO space  $\widehat{V}_{h,0}^{k+}$ reads
$$\widehat{V}_{h}^{k+}=\{((v_{T})_{T \in \mathcal{T}_h},(v_{F})_{F \in \mathcal{F}_h}) \mid v_{T} \in \mathbb{P}_{d}^{k+1}(T), v_{F} \in \mathbb{P}_{d-1}^{k}(F), k\geq 0\}.$$ 
Impose the Dirichlet boundary condition to get  
$\widehat{V}_{h,0}^{k+}=\{\widehat{w}_h\in\widehat{V}_{h}^{k+}\mid v_F\equiv 0, \text{for}\; F\in \mathcal{F}^b \}.$ 
\subsection{Norms}\label{sec_23}
We denote the norm (resp. inner product) on the space $L^2(\Omega)$ as $\norm{\bullet}$ (resp. $(\bullet,\bullet)$ ). The norm (resp. inner product) on the space $L^2(T)$ and $L^2(F)$ as $\norm{\bullet}_T$ and $\norm{\bullet}_F$  (resp. $(\bullet,\bullet)_T$ and $(\bullet,\bullet)_F$).  
For all $\widehat{v}_{T}:=(v_{T}, (v_F)_{F\in \mathcal{F}_T}) \in \widehat{V}_{T}^{k},$ define a  $H^{1}$-like semi-norm $\widehat{V}_{T}^{k}$ as 
\begin{align*}
|\widehat{v}_{T}|_{1,T}^{2}:=\norm{\nabla v_{T}}_{T}^{2}+\sum_{F\in \mathcal{F}_T}h_{T}^{-1}\norm{v_{T}-v_{F}}_{F}^{2}. 
\end{align*}
The norm on the global HHO space $\widehat{V}_{h, 0}^{k}$ is defined as 
\begin{align*}
      \hhonorm{\hhoelem{w}}^2 :=\sum_{T\in \mathcal{T}_h}(\lVert\nabla v_{T}\rVert_{T}^{2}+h_{T}^{-1}\sum_{F\in \mathcal{F}_T}\lVert v_{T}-v_{F}\rVert_{F}^{2})\quad \text{for all}\; \hhoelem{w}\in \widehat{V}_{h, 0}^{k}.
\end{align*}
\subsection{Local and global reduction operators}
For any given non-negative integer $k$ and given $T\in{\cal T}_{h}$(resp. $F\in {\cal F}_{h}$) the {\it local $L^{2}$-projection} $\Pi_{T}^{k}: L^{2}(T)\rightarrow \mathbb{P}_{d}^{k}(T)$(resp. $\Pi_{F}^{k}: L^{2}(F)\rightarrow \mathbb{P}_{d-1}^{k}(F)$) reads: For given $v\in L^{2}(T)$(resp. $v\in L^{2}(T)$),
\begin{equation}
  (\Pi_{T}^{k}v,w)_{T}=(v,w)_{T}\; \forall w\in \mathbb{P}_{d}^{k}(T) \quad \text{and} \quad (\Pi_{F}^{k}v,w)_{F}=(v,w)_{F}\; \forall w\in \mathbb{P}_{d-1}^{k}(F).\label{L_2_proj_map}
\end{equation}
The {\it local reduction operator} $\widehat{I}_{T}^{k}: H^{1}(T) \rightarrow\widehat{V}_{T}^{k}$ is defined as: For all $v \in H^{1}(T)$,
\begin{equation*}
\widehat{I}_{T}^{k} v:=(\Pi_{T}^{k} v,(\Pi_{F}^{k} v)_{F \in \mathcal{F}_{T}}).
\end{equation*}
For all $v \in H^{1}(\Omega)$, Define {\it global reduction operator} $\widehat{I}_{h}^{k}: H^{1}(\Omega) \rightarrow\widehat{V}_{h}^{k}$ as
$$
\widehat{I}_{h}^{k} v:=((\Pi_{T}^{k} v)_{T \in \mathcal{T}_{h}},(\Pi_{F}^{k} v)_{F \in \mathcal{F}_{h}}).
$$
\subsection{Local reconstruction and stabilization operators.} 
The {\it local reconstruction operator} $\mathcal{R}_{T}: \widehat{V}_{T}^{k} \rightarrow \mathbb{P}_{d}^{k+1}(T)$ is defined as: For any given $\widehat{v}_{T}=(v_{T},(v_{F})_{F \in \mathcal{F}_{T}})$ and $T \in \mathcal{T}_{h},$
\begin{align}
(\nabla \mathcal{R}_{T} \widehat{v}_{T}, \nabla w)_{T} & =(\nabla v_{T}, \nabla w)_{T}+\sum_{F \in \mathcal{F}_{T}}(v_{F}-v_{T}, \nabla w \cdot \boldsymbol{n}_{T F})_{F} \quad \forall w \in \mathbb{P}_{d}^{k+1}(T),\label{eqn:recon1} \\
(\mathcal{R}_{T} \widehat{v}_{T}, 1)_{T} & =(v_{T}, 1)_{T}.\label{eqn:recon2}
\end{align}
 Define a {\it global reconstruction operator} $\mathcal{R}_{h}$ : $\widehat{V}_{h}^{k} \rightarrow\mathbb{P}_{d}^{k+1}(\mathcal{T}_{h})$ by $\mathcal{R}_{h} \widehat{v}_{h}|_{T}=\mathcal{R}_{T} \widehat{v}_{T}$. The  {\it local stabilization operator} $S_F:\widehat{V}_T^k\rightarrow\mathbb{P}_{d-1}^k(F)$ reads $S_{F}(\lochhoelem{v}) = \Pi^k_{F}(v_T|_F-v_F+((I-\Pi^k_T)\mathcal{\mathcal{R}}_T \lochhoelem{v})|_{F}).$
\subsection{Local elliptic projection operator}
Denote $\mathcal{E}_{T}^{k+1}=\mathcal{R}_{T} \circ \widehat{I}_{K}^{k}$ where $\mathcal{E}_{T}^{k+1}:H^{1}(T)\rightarrow \mathbb{P}_{d}^{k+1}(T)$ is the elliptic operator uniquely defined such that for all $v\in H^{1}(T),$
\begin{align}
&(\nabla\mathcal{E}_{T}^{k+1}(v),\nabla q)_{T}=(\nabla(v),\nabla q)_{T}\quad \forall q\in \mathbb{P}_{d}^{k+1}(T)^{\perp},\label{elli_proj_1}\\
&(\mathcal{E}_{T}^{k+1}(v),1)_{T}=(v,1)_{T}\label{elli_proj_1}.
\end{align}
Further, for $v\in H^{k+1}(T),$ it holds
\begin{equation}
    \norm{\nabla(v-\mathcal{E}_{T}^{k+1}(v))}_{T}\leq  \norm{\nabla(v-\Pi_{T}^{k+1}(v))}_{T}\leq Ch^{k}|v|_{k+1}.\label{ellip_prj_bd}
\end{equation}

\subsection{Bilinear forms}
For all $v,w \in H^{1}_{0}(\Omega),$ define $a(v, w) =(\nabla v,\nabla w).$ The discrete {\it bilinear form} $a_h$ is defined as: For all $\widehat{v}_h, \widehat{w}_h\in \widehat{V}_{h,0}^{k},$
\begin{align*}
  a_h(\widehat{v}_h, \widehat{w}_h) : = \sum_{T\in \mathcal{T}_h}( ( \nabla \mathcal{R}_T\widehat{v}_T, \nabla \mathcal{R}_T\widehat{w}_T)_T + S_T(\widehat{y}_T, \widehat{w}_T)),  
\end{align*}
 where local stabilization map $S_T$ is defined as
 \begin{align*}
S_T(\lochhoelem{y},\lochhoelem{w})=h_T^{-1}\sum_{F\in \mathcal{F}_T}(S_{F}(\lochhoelem{y}), S_{F}(\lochhoelem{w}))_{F}.
 \end{align*}

\subsection{Trace inequality and estimate for projection}
 There exist real numbers $C_{\mathrm{tr}}$ and $C_{\mathrm{tr}, \mathrm{c}}$ depending on $\varrho$ but independent of mesh parameter $h$ such that the following discrete and continuous trace inequalities hold for all $T \in \mathcal{T}_{h}$ and $F \in \mathcal{F}_{T}$ (see, \cite[Lemmas 1.46 and 1.49]{Ern_DiPietro_Book2012}):
\begin{align}
\|v\|_{F} & \leq C_{\mathrm{tr}} h_{T}^{-1 / 2}\|v\|_{T} \quad \forall v \in \mathbb{P}_{d}^{l}(T)\label{Inv_ineq} \\
\|v\|_{\partial T} & \leq C_{\mathrm{tr}, \mathrm{c}}(h_{T}^{-1}\|v\|_{T}^{2}+h_{T}\|\nabla v\|_{T}^{2})^{1 / 2} \quad \forall v \in H^{1}(T),\label{trace_ineq}
 \end{align}
where $\mathbb{P}_{d}^{l}(T)$ is the space of polynomial of degree at most $l$ on $T \in \mathcal{T}_{h}$. There exists a real number $C_{\text {app }}$ depending on $\varrho$ and $l$ but independent of $h$ such that, for all $T \in \mathcal{T}_{h}$, denoting by $\Pi_{T}^{l}$ the $L^{2}$-orthogonal projector on $\mathbb{P}_{d}^{l}(T)$, the following holds (see \cite[Lemmas 1.58 and 1.59]{Ern_DiPietro_Book2012}): For all $s \in\{1, \ldots, l+1\}$ and all $v \in H^{s}(T)$,
\begin{equation}
|v-\Pi_{T}^{l} v|_{H^{m}(T)}+h_{T}^{1 / 2}|v-\Pi_{T}^{l} v|_{H^{m}(\partial T)} \leq C_{\mathrm{app}} h_{T}^{s-m}|v|_{H^{s}(T)} \quad \forall m \in\{0, \ldots, s-1\} ,\label{app_prop_1}
\end{equation}
where $|\cdot|_{H^{m}(\partial T)}$ denotes the face-wise $H^{m}$ semi-norm when the boundary $\partial T$ of an element $T \in \mathcal{T}_{h}$ is written as a union of faces.

\begin{remark}
The norm, {\it local/global reduction}, {\it local/global reconstruction}, and {\it local/global stabilization operators} for the space $\widehat{V}_{h}^{k+}$ are defined analogously to $\widehat{V}_{h}^{k}.$ 
\end{remark}
\subsection{Integration by parts formula}
For any $T\in {\cal T}_{h}$ and sufficiently smooth functions $v$ and $w$, an {\it integration by parts} leads to  
\begin{equation}
(\Delta w,v)_{T}=-(\nabla w,\nabla v)_{T}+\sum_{F\in \partial T}(v,\partial_{\bf n}w)_{F}.\label{int_part_2}
\end{equation}
\section{Optimal control problem without control constraints}\label{sec:OCP_HHO_UC}
\subsection{Existence, uniqueness, and regularity for \eqref{eqn:conti_costfun}-\eqref{eqn:conti_state}}
In this subsection, the unconstrained distributed optimal control problem \eqref{eqn:conti_costfun}-\eqref{eqn:conti_state} is considered with control in  $L^2(\Omega).$ Existence and uniqueness of the solution are standard and can be found in \cite[Theorem 2.14]{trolzstch:2005:Book}.
The problem is formulated as the following optimality system using first-order necessary optimality condition: Seek $(y, u, \phi)\in H^1_0(\Omega) \times  L^2(\Omega)\times H^1_0(\Omega) $, such that
    \begin{align}
         a(y, w)&=(f, w)+(u, w) \quad \forall w \in H_{0}^1(\Omega), \label{eqn:conti_os_state}\\ a(w, \phi)&=(y-y_d, w)  \quad\forall w \in H^1_{0}(\Omega),\label{eqn:conti_os_adj_state}\\ 
        u &= -\phi/\lambda. \label{eqn:conti_os_control}  
    \end{align}
\begin{lemma}[Regularity]\label{lem:regularity_UC}
For $f,y_d\in H^k(\mathcal{T}_h)$ with $k\geq 0,$ the solution of  \eqref{eqn:conti_os_state}-\eqref{eqn:conti_os_control} has the regularity $y,u, \phi \in H^{k+2}(\mathcal{T}_h).$
\end{lemma}
\begin{proof}
The source term $f\in H^k(\mathcal{T}_h)$ and $u\in L^2(\Omega),$ thus $f+u\in L^2(\Omega),$ hence from the equation \eqref{eqn:conti_state} we have $y\in H^2(\Omega).$ Given that $y_d\in H^k(\mathcal{T}_h)$, thus $y-y_d\in H^2(\mathcal{T}_h),$ therefore from the strong form of the adjoint state equation \eqref{eqn:conti_os_adj_state}, we have $\phi\in H^4(\mathcal{T}_h).$ Thus equation \eqref{eqn:conti_os_control}, yields $u\in H^4(\mathcal{T}_h).$ Now $f+u\in H^4(\mathcal{T}_h),$ hence $y\in H^6(\mathcal{T}_h).$ Then $\phi \in H^8(\mathcal{T}_h)$ and hence $u \in H^8(\mathcal{T}_h).$ Thus, we keep on doing this bootstrap argument to get solutions $y,u, \phi \in H^{k+2}(\mathcal{T}_h),$ and this proves the lemma.
\end{proof}
\subsection{Scheme UC1}\label{Scheme:UC_1}
Introduce the discrete control space $U_h :=\{ v_h\in L^2(\Omega)| \; v_h|_T \in \mathbb{P}_k(T) \; \forall T \in \mathcal{T}_h\}.$  The {\it HHO scheme} for \eqref{eqn:conti_costfun}-\eqref{eqn:conti_state} is given by
 \begin{equation}\label{eqn:disc_costfun_UM1}
 \min_{\widehat{y}_h \in \Vsp, u_h\in U_h} J(\widehat{y}_h, u_h):=\frac{1}{2}\|y_{\mathcal{T}}-y_{d}\|^{2}+\frac{\lambda}{2}\|u_h\|^{2}, 
\end{equation} 
subject to the PDE, find $\hhoelem{y}\in \widehat{V}_{h,0}^{k} $ such that
\begin{equation}\label{eqn:disc_state_UM1}
 a_h(\widehat{y}_h, \widehat{w}_h)=(f,w_{\mathcal{T}})+(u_h, w_{\mathcal{T}}) \quad \forall \widehat{w}_h \in \widehat{V}_{h,0}^{k}.   
\end{equation}
Next, the {\it reduced cost functional} is introduced, which will be useful for the existence and uniqueness of the discrete solution and to derive the optimality system.

\medskip\noindent
{\it Reduced cost functional}: Define $\widehat{\mathfrak{S}}_h: L^2(\Omega) \times U_h \rightarrow\Vsp$ such that $\widehat{\mathfrak{S}}_h (f,u_h) = \hhoelem{y}$ be the solution operator of \eqref{eqn:disc_state_UM1}. The {\it reduced cost functional} $j_h$ reads
\begin{align}\label{eqn:reduced_costfun_UM1}
   j_h(u_h)=\frac{1}{2}\|\mathfrak{S}_{\mathcal{T}}(f,u_h)-y_{d}\|^{2}+\frac{\lambda}{2}\|u_h\|^{2}, 
\end{align}
where $\mathfrak{S}_{\mathcal{T}}(f,u_h)$ be the volume part of $\widehat{\mathfrak{S}}_h (f,u_h).$ The discrete problem \eqref{eqn:disc_costfun_UM1}-\eqref{eqn:disc_state_UM1} with $\widehat{\mathfrak{S}}_h$ leads to
\begin{align}\label{eqn:reduced_prob_UM1}
  \min_{u_h \in U_h} j_h(u_h)=\frac{1}{2}\|\mathfrak{S}_{\mathcal{T}}(f,u_h)-y_{d}\|^{2}+\frac{\lambda}{2}\|u_h\|^{2}.
\end{align}
The following theorem establishes the existence and uniqueness of the discrete problem \eqref{eqn:reduced_prob_UM1}.
\begin{theorem}[Existence and uniqueness of discrete solution]\label{thm:existance_thm}
There exists a unique solution of the problem \eqref{eqn:reduced_prob_UM1}.
\end{theorem}

\begin{proof}
The solution operator $\widehat{\mathfrak{S}}_h: U_h \rightarrow\widehat{V}^k_{h,0}$ of the PDE \eqref{eqn:disc_state_UM1} is continuous i.e., $\hhonorm{\hhoelem{y}}=\hhonorm{\widehat{\mathfrak{S}}_h (f,u_h)} \lesssim \norm{f}+ \norm{u_h},$ using this stability estimate of the PDE \eqref{eqn:disc_state_UM1}. Using Poincar\'e inequality we have $\norm{y_{\mathcal{T}}}\leq C_p \hhonorm{\hat{y}_h}$. Therefore $\widehat{\mathfrak{S}}_h: L^2(\Omega)\times U_h \rightarrow L^2(\Omega)$ is a continuous affine linear operator. The norm  $\norm{\cdot}$ is continuous and convex, therefore $j_h$ is continuous and convex. Clearly, $j_h(u_h)\geq 0,$ therefore the infimum exists. Let $\alpha = \inf j_h(u_h).$ Let $u_h^n\in U_h$ be an infimizing sequence such that $j_h(u^n_h)$ converges to $\alpha$ as $n$ tends to infinity. Thus the sequence $u_h^n$ is bounded in $L^2.$ Then there exists a subsequence $u_h^n$ still indexed by $n$ such that it is converges to $u_h^*$ weakly in $L^2(\Omega)$. It is clear that $u_h^* \in U_h$ since $U_h$ is weakly closed. The discrete cost functional $j_h$ is continuous and convex hence it is weakly lower semicontinuous. Thus $j_h(u_h^*) \leq \liminf j_h(u_h^n)$ and hence $j_h(u_h^*) \leq \alpha.$ Therefore $u_h^*$ be a minimizer. The uniqueness follows from the strict convexity of the reduced cost functional $j_h.$   
\end{proof}

\noindent
The discrete optimality system reads as: find $(\widehat{y}_h, u_h, \widehat{\phi}_h)\in \widehat{V}_{h,0}^{k} \times U_h \times \widehat{V}_{h,0}^{k} $ such that
\begin{align}
     a_h(\widehat{y}_h, \widehat{w}_h)&=(f, w_{\mathcal{T}})+(u_h, w_{\mathcal{T}}) \quad \forall \widehat{w}_h \in \widehat{V}_{h,0}^k,\label{eqn:disc_state}\\
   a_h(\widehat{w}_h, \widehat{\phi}_h)&=(y_{\mathcal{T}}, w_{\mathcal{T}})-(y_d, w_{\mathcal{T}} ) \quad \forall \widehat{w}_h \in \widehat{V}_{h, 0}^k,\label{eqn:disc_adjstate}\\
   u_h &=-\frac{1}{\lambda} \phi_{\mathcal{T}}. \label{eqn:disc_VI} 
\end{align}
The system \eqref{eqn:disc_state}-\eqref{eqn:disc_VI} can be derived from the Lagrangian functional defined as
\begin{align}\label{eqn:disc_Lagrangian}
    \mathcal{L}_h(\widehat{y}_h, u_h, \widehat{\phi}_h)= \frac{1}{2}\|y_{\mathcal{T}}-y_{d}\|^2+\frac{\lambda}{2}\|u_h\|^2-a_h(\widehat{y}_h, \widehat{\phi}_h) +(f, \phi_{\mathcal{T}})+(u_h, \phi_{\mathcal{T}}).
\end{align}
Taking derivative of $\mathcal{L}_h$ with respect to $\hhoelem{y}$ (resp. $u_h$) and equating to zero i.e., 
$\mathcal{D}_{\hhoelem{y}}\mathcal{L}_h(\widehat{y}_h, u_h, \widehat{\phi}_h)\hhoelem{w} = 0 \; \forall \hhoelem{w} \in \Vsp
$ (resp. $\mathcal{D}_{u_h}\mathcal{L}_h(\widehat{y}_h, u_h, \widehat{\phi}_h)w_h = 0 \;\forall w_h \in U_h$) yields \eqref{eqn:disc_adjstate} (resp. \eqref{eqn:disc_VI}).

\medskip\noindent
To derive the error estimates, we require the following auxiliary problems and second-order optimality condition.
\begin{itemize}
    \item (\it Auxiliary problems). Seek $\tilde{y}_h(u),\tilde{\phi}_h(u) \in \widehat{V}_{h, 0}^{k}$ such that
\begin{align}
 a_h(\tilde{y}_h(u), \widehat{w}_h)&=(f, w_{\mathcal{T}} )+(u, w_{\mathcal{T}} ) \quad \forall \widehat{w}_h \in \widehat{V}_{h, 0}^{k},\label{eqn:aux_y(u)_UM1} \\
a_h(\widehat{w}_h, \tilde{\phi}_h(u))&=(y-y_d, w_{\mathcal{T}} ) \quad \forall \widehat{w}_h \in \widehat{V}_{h, 0}^{k}\label{eqn:aux_p(u)_UM1}.
\end{align}
\item (\it Second-order optimality condition).  For all $v_h \in U_h,$ 
 \begin{align}\label{eqn:disc_red_costfun_UM1}
 j_h'(u_h)v_h &= (\mathfrak{S}_{\mathcal{T}}(f,u_h)-y_{d},\mathfrak{S}_{\mathcal{T}}(0,v_h))+\lambda(u_h,v_h),\nonumber\\
 j_h''(u_h)v_h^2 &= \norm{\mathfrak{S}_{\mathcal{T}}(0,v_h)}^2+\lambda \norm{v_h}^2.
\end{align} 
From \eqref{eqn:disc_red_costfun_UM1}, it holds
\begin{align}\label{eqn:SSC}
    j_h''(u_h)v_h^2 \geq \lambda \norm{v_h}^2.
\end{align}
\end{itemize}
\begin{theorem}[Abstract $L^2$- error estimates of control]\label{thm:L2_control_error}
Let $u$ (resp. $u_h$) solves \eqref{eqn:conti_costfun}-\eqref{eqn:conti_state} (resp. \eqref{eqn:reduced_prob_UM1}). There holds,
\begin{align*}
     \norm{u-u_h} \leq &C(C_p,\lambda,\alpha_0) \big(\norm{u - v_h} + \norm{y - \tilde{y}_{\mathcal{T}} (u)}\\&+ \norm{\phi -\tilde{\phi}_{\mathcal{T}} (u)}\big)\quad \forall v_h \in U_h,   
\end{align*}
where $C(C_p,\lambda,\alpha_0)$ be the constant depends on the Poincar\'e constant $C_p,$ regularization parameter $\lambda,$ and ellipticity constant $\alpha_0.$
\end{theorem}
\begin{proof}
The triangle inequality yields
    \begin{align}\label{tria_ineq}
        \norm{u-u_h} \leq \norm{u-v_h} +\norm{v_h-u_h},
    \end{align}
    for all $v_h \in U_h,$
 Now the second term  $\norm{v_h-u_h}$ on the right-hand side of \eqref{tria_ineq} is controlled below. 

\medskip\noindent
A standard computation produces  
$j_h'(v_h)(v_h-u_h)=(\phi_{\mathcal{T}}(v_h)+\lambda v_h, v_h-u_h)$ where $\widehat{\phi}_h(v_h)$ solves 
\begin{align}\label{eqn:auxproblem3}
    a_h(\widehat{w}_h,\widehat{\phi}_h(v_h)) = (y_{\mathcal{T}}(v_h)-y_d, w_{\mathcal{T}}) \quad \forall \widehat{w}_h \in \widehat{V}^k_{h,0}.
\end{align}
The second-order sufficient condition \eqref{eqn:SSC}, mean-value property of $j_h'(\bullet),$ \eqref{eqn:disc_VI}, the quantity $j_h'(v_h)(v_h-u_h)$ computed above, and elementary algebra lead to
\begin{align*}
    \lambda \norm{v_h-u_h}^2 &\leq j_h''(w_h)(v_h- u_h)^2
     = j_h'(v_h)(v_h-u_h)- j_h'(u_h)(v_h-u_h)
    = (\phi_{\mathcal{T}}(v_h)+\lambda v_h, v_h-u_h)\nonumber\\&
    = (\phi_{\mathcal{T}}(v_h)-\phi, v_h-u_h)+\lambda(v_h-u,v_h-u_h)+(\phi+\lambda u, v_h-u_h).
     \end{align*}
 Note that $(\phi+\lambda u, v_h-u_h)=0$ from first-order necessary optimality condition \eqref{eqn:conti_os_control}  and 
Cauchy-Schwarz inequality reveal
\begin{align}\label{auxestm:03}
    \lambda\norm{v_h-u_h} \leq \norm{\phi-\phi_{\mathcal{T}}(v_h)} + \lambda \norm{u - v_h}.
\end{align}
Now we estimate $\norm{\phi-\phi_{\mathcal{T}}(v_h)}$ as follows.
The use of triangle inequality yields 
\begin{align}\label{auxestm:04}
   \norm{\phi-\phi_{\mathcal{T}}(v_h)} \leq \norm{\phi - \tilde{\phi}_{\mathcal{T}}(u)} + \norm{\tilde{\phi}_{\mathcal{T}}(u)-\phi_{\mathcal{T}}(v_h)}. 
\end{align}
A subtraction of \eqref{eqn:auxproblem3} from \eqref{eqn:aux_p(u)_UM1} leads to
    $a_h(\widehat{w}_h, \tilde{\phi}_h(u)-\widehat{\phi}_h(v_h)) = (y-y_{\mathcal{T}}(v_h),w_{\mathcal{T}}) \;  \forall \widehat{w}_h \in \widehat{V}^k_{h,0}.$
The choice $\widehat{w}_h = \tilde{\phi}_h(u)-\widehat{\phi}_h(v_h)$ and ellipticity of $a_{h}(\bullet,\bullet)$ in the left-hand side and Cauchy-Schwarz inequality in the right-hand side of this expression show
    $\alpha_0 \norm{\tilde{\phi}_h(u)-\widehat{\phi}_h(v_h)}^2_{1,h} \leq \norm{y-y_{\mathcal{T}}(v_h)} \norm{\tilde{\phi}_{\mathcal{T}}(u)-\widehat{\phi}_{\mathcal{T}}(v_h)}.$
Apply Poincar\'e inequality 
to get
\begin{align}\label{auxestm:05}
    \alpha_0 \norm{\tilde{\phi}_{\mathcal{T}}(u)-\widehat{\phi}_{\mathcal{T}}(v_h)} \leq C_p \norm{y-y_{\mathcal{T}}(v_h)}.
\end{align}
Introduce the split to get 
  $  \norm{y-y_{\mathcal{T}}(v_h)} \leq \norm{y-\tilde{y}_{\mathcal{T}}(u)} + \norm{\tilde{y}_{\mathcal{T}}(u)-y_{\mathcal{T}}(v_h)}.$
Recall that $\widehat{y}_h(v_h)$ solves
    $a_h(\widehat{y}_h(v_h),\widehat{w}_h) = (f, w_{\mathcal{T}}) + (v_h, w_{\mathcal{T}}) \;\forall \widehat{w}_h \in \widehat{V}^k_{h,0}.$
These two results reveal $C_p^2 \alpha_0 \norm{\tilde{y}_{\mathcal{T}}(u)-y_{\mathcal{T}}(v_h)}\leq \norm{u-v_h}.$ 
This with \eqref{auxestm:05} and elementary algebra lead to
\begin{align}\label{auxestm:25}
     \norm{\tilde{\phi}_{\mathcal{T}}(u)-\widehat{\phi}_{\mathcal{T}}(v_h)} \leq C(C_p,\alpha_0) (\norm{y-\tilde{y}_{\mathcal{T}}(u)} + \norm{u-v_h}).
\end{align}
A combination of \eqref{auxestm:04} and \eqref{auxestm:25} establish
\begin{align}\label{auxestm:24}
   \norm{\phi-\phi_{\mathcal{T}}(v_h)} \leq C(C_p,\alpha_0) (\norm{\phi - \tilde{\phi}_{\mathcal{T}}(u)} + \norm{y-\tilde{y}_{\mathcal{T}}(u)} + \norm{u-v_h}). 
\end{align}
Utilize \eqref{auxestm:24} in \eqref{auxestm:03}, and the final result in \eqref{tria_ineq} concludes the lemma.
\end{proof}
\begin{lemma}\label{lem:estm_aux_state}
Under the regularity of solution in Lemma \ref{lem:regularity_UC}, it holds
   $$\norm{y-\tilde{y}_{\mathcal{T}}(u)}=\mathcal{O}(h^{k+1}).$$ 
\end{lemma}
\begin{proof}
Introduce the split to get
  \begin{align}\label{estm007_15}
      \norm{y - \tilde{y}_{\mathcal{T}} (u)} \leq \norm{y-\Pi^k_{\mathcal{T}} y} + \norm{\Pi^k_{\mathcal{T}} y - \tilde{y}_{\mathcal{T}} (u)}.
  \end{align}
Now, we control both the terms on the right-hand side of the above inequality. The proof is divided into two steps. 

\noindent\medskip  
\noindent \medskip{\it Step 1.} (Control for $\norm{y-\Pi^k_{\mathcal{T}} y}$) The bound of first term  
    $\norm{y-\Pi^k_{\mathcal{T}} y}\lesssim h^{k+1}\norm{y}_{H^{k+2}(\mathcal{T}_h)}$ 
 follows from \eqref{app_prop_1}.
  
\noindent \medskip{\it Step 2.} (Control for $\norm{\Pi^k_{\mathcal{T}} y - \tilde{y}_{\mathcal{T}} (u)}$) 
 Let $\sigma_{\mathcal{T}}:=\Pi^k_{\mathcal{T}} y - \tilde{y}_{\mathcal{T}} (u).$
For any $\Psi_{\sigma} \in H_{0}^{1}(\Omega)$ define
\begin{align}\label{aux:eqn007_13}
     a(v, \Psi_{\sigma})=( \sigma_{\mathcal{T}},v) \quad\; \forall \; v \in H_{0}^{1}(\Omega).
\end{align}
Note that $-\Delta \Psi_{\sigma}=\sigma_{\mathcal{T}}.$ An application of integration by parts formula \eqref{int_part_2} leads to
\begin{align}\label{estm_007_13}
    \norm{\sigma_{\mathcal{T}}}^2=-(\sigma_{\mathcal{T}}, \Delta \Psi_{\sigma})=\sum_{T \in \mathcal{T}}((\nabla \sigma_{T}, \nabla \Psi_{\sigma})_{T}+(\sigma_{\partial T}-\sigma_{T}, \boldsymbol{n}_{T} \cdot \nabla \Psi_{\sigma})_{\partial T}),
\end{align}
where we used that $\Psi_{\sigma} \in H^{1+s}(\Omega), s>\frac{1}{2}$ and $\partial_{{\bf n}}\Psi_{\sigma}= \boldsymbol{n}_{T} \cdot \nabla \Psi_{\sigma}$ is single-valued across the mesh interface and  $\sigma_{F}=0$ for all $F \in \mathcal{F}^{b}.$

\noindent \medskip Denote $\xi:=\Psi_{\sigma}-\mathcal{E}_{\mathcal{T}}^{k+1}(\Psi_{\sigma})$. Add and subtract $\mathcal{R}_{T}(\widehat{I}_{T}^{k}(\Psi_{\sigma}))=\mathcal{E}_{T}^{k+1}(\Psi_{\sigma})$ for all $T\in\mathcal{T}$ in the right-hand side of \eqref{estm_007_13},  definition of $\mathcal{R}_{T}(\widehat{\sigma}_{T}),$ and $(\nabla \sigma_{T}, \nabla \xi)_{T}=0$ reveal
\begin{align}
    \|\sigma_{\mathcal{T}}\|^{2}&=\sum_{T \in \mathcal{T}}(\sigma_{\partial T}-\sigma_{T}, \boldsymbol{n}_{T} \cdot \nabla \xi)_{L^{2}(\partial T)}+\sum_{T \in \mathcal{T}}(\nabla \mathcal{R}_{T}(\widehat{\sigma}_{T}), \nabla \mathcal{R}_{T}(\widehat{I}_{T}^{k}(\Psi_{\sigma})))_{T}.\nonumber
    \end{align}
The definitions of $a_{h}(\bullet,\bullet),$ $\widehat{\sigma}_{h},$ and  \eqref{eqn:aux_y(u)_UM1} lead to
\begin{align}
    \|\sigma_{\mathcal{T}}\|^{2}&= \sum_{T \in \mathcal{T}}(\sigma_{\partial T}-\sigma_{T}, \boldsymbol{n}_{T} \cdot \nabla \xi)_{L^{2}(\partial T)}-s_{h}(\widehat{\sigma}_{h}, \widehat{I}_{h}^{k}(\Psi_{\sigma}))+a_{h}(\widehat{\sigma}_{h}, \widehat{I}_{h}^{k}(\Psi_{\sigma}))\nonumber\\
    &= \sum_{T \in \mathcal{T}}(\sigma_{\partial T}-\sigma_{T}, \boldsymbol{n}_{T} \cdot \nabla \xi)_{L^{2}(\partial T)}-s_{h}(\widehat{\sigma}_{h}, \widehat{I}_{h}^{k}(\Psi_{\sigma}))+(f+u, \Pi_{\mathcal{T}}^{k}(\Psi_{\sigma}))-a_{h}(\widehat{I}_{h}^{k}(y), \widehat{I}_{h}^{k}(\Psi_{\sigma})).\nonumber
    \end{align}
An application of \eqref{eqn:conti_state} in above equation, definition of $\Pi_{\mathcal{T}}^{k}$, and elementary algebra establish   
    \begin{align}\label{estm:sigma}
   \|\sigma_{\mathcal{T}}\|^{2}= &\big(\sum_{T \in \mathcal{T}}(\sigma_{\partial T}-\sigma_{T}, \boldsymbol{n}_{T} \cdot \nabla \xi)_{L^{2}(\partial T)}-s_{h}(\widehat{\sigma}_{h}, \widehat{I}_{h}^{k}(\Psi_{\sigma}))\big)+\big((\nabla y, \nabla \Psi_{\sigma})-a_{h}(\widehat{I}_{h}^{k}(y), \widehat{I}_{h}^{k}(\Psi_{\sigma}))\big)\nonumber\\&-\big(((f+u)-\Pi_{\mathcal{T}}^{k}(f+u), \Psi_{\sigma}-\Pi_{\mathcal{T}}^{k}(\Psi_{\sigma})\big)
    :=\Xi_{1}+\Xi_{2}-\Xi_{3}.
\end{align}
It remains to bound $\Xi_{i},~i=1,2,3$. The Cauchy-Schwarz inequality, Lemma \ref{Lem:auxilary02} (to bound $s_{h}(\widehat{\sigma}_{h}, \widehat{\sigma}_{h})^{\frac{1}{2}}$, and Lemma \ref{Lem:stabilizer_estm} (to bound $s_{h}(\widehat{I}_{h}^{k}(\Psi_{\sigma}), \widehat{I}_{h}^{k}(\Psi_{\sigma}))^{\frac{1}{2}}$ give
$$
|\Xi_{1}| \leq C\|\widehat{\sigma}_{h}\|_{\widehat{V}_{h, 0}^{k}}( \sum_{T\in \mathcal{T}_h} (\norm{\nabla(\Psi_{\sigma}-\mathcal{E}_{\mathcal{T}}^{k+1}(\Psi_{\sigma}))}_T^2 + h_T\norm{\nabla(\Psi_{\sigma}-\mathcal{E}_{\mathcal{T}}^{k+1}(\Psi_{\sigma}))}_{\partial T}^2))^{1/2}.
$$   
The approximation property of the elliptic projection gives 
\begin{align*}
    ( \sum_{T\in \mathcal{T}_h} (\norm{\nabla(\Psi_{\sigma}-\mathcal{E}_{\mathcal{T}}^{k+1}(\Psi_{\sigma}))}_T^2 + h_T\norm{\nabla(\Psi_{\sigma}-\mathcal{E}_{\mathcal{T}}^{k+1}(\Psi_{\sigma}))}_{\partial T}^2))^{1/2}\leq C h^{s}|\Psi_{\sigma}|_{H^{1+s}(\Omega)}.
\end{align*}
Note that $\|\Psi_{\sigma}\|_{H^{1+s}(\Omega)} \leq C_{\text {ell }} \|\sigma_{\mathcal{T}}\|$ by the elliptic regularity result from \eqref{aux:eqn007_13}. The Lemma \ref{energy_error_estimate} to bound $\|\widehat{\sigma}_{h}\|_{\widehat{V}_{h, 0}^{k}}$, we infer that
\begin{align*}
|\Xi_{1}| &\leq C  h^{s}\|\sigma_{\mathcal{T}}\|( \sum_{T\in \mathcal{T}_h} (\norm{\nabla(y-\mathcal{E}_{\mathcal{T}}^{k+1}y)}_T^2 + h_T\norm{\nabla(y-\mathcal{E}_{\mathcal{T}}^{k+1}y)}_{\partial T}^2))^{1/2}.
\end{align*}
Next apply the definitions of $a_{h}(\bullet,\bullet)$, $\mathcal{R}_{\mathcal{T}} \circ \widehat{I}_{h}^{k}=\mathcal{E}_{\mathcal{T}}^{k+1}$, and the orthogonality property \eqref{elli_proj_1} to get
$$
\begin{aligned}
\Xi_{2} & =(\nabla y, \nabla \zeta_e)_{\boldsymbol{L}^{2}(\Omega)}-(\nabla_{\mathcal{T}} \mathcal{E}_{\mathcal{T}}^{k+1}(y), \boldsymbol{\nabla}_{\mathcal{T}} \mathcal{E}_{\mathcal{T}}^{k+1}(\Psi_{\sigma}))-s_{h}(\widehat{I}_{h}^{k}(y), \widehat{I}_{h}^{k}(\Psi_{\sigma})), \\
& =(\nabla_{\mathcal{T}}(y-\mathcal{E}_{\mathcal{T}}^{k+1}(y)), \nabla_{\mathcal{T}}(\Psi_{\sigma}-\mathcal{E}_{\mathcal{T}}^{k+1}(\Psi_{\sigma})))_{\boldsymbol{L}^{2}(\Omega)}-s_{h}(\widehat{I}_{h}^{k}(y), \widehat{I}_{h}^{k}(\Psi_{\sigma})).
\end{aligned}
$$
The Cauchy-Schwarz inequality and Lemma  \ref{Lem:stabilizer_estm} reveal
$|\Xi_{2}| \leq C\|\nabla_{\mathcal{T}}(y-\mathcal{E}_{\mathcal{T}}^{k+1}(y))\|\|\nabla_{\mathcal{T}}(\Psi_{\sigma}-\mathcal{E}_{\mathcal{T}}^{k+1}(\Psi_{\sigma}))\|.$
An application of \eqref{ellip_prj_bd} and \eqref{aux:eqn007_13} establish
$|\Xi_{2}| \leq C h^{s} \|\nabla_{\mathcal{T}}(y-\mathcal{E}_{\mathcal{T}}^{k+1}(y))\|\|\sigma_{\mathcal{T}}\|.$

\medskip\noindent
The definition of $\Pi_{\mathcal{T}}^{k}$ leads to
 $\Xi_{3} \leq \lVert (f+u)-\Pi_{\mathcal{T}}^{k}(f+u)\rVert \lVert \Psi_{\sigma}-\Pi_{\mathcal{T}}^{k}\Psi_{\sigma}\rVert
\leq C \lVert (f+u)-\Pi_{\mathcal{T}}^{k}(f+u)\rVert  h_{T}^{1+s}\norm{\Psi_{\sigma}}_{H^{1+s}(T)}
\leq C h_{T}^{1+s} \lVert (f+u)-\Pi_{\mathcal{T}}^{k}(f+u)\rVert  \norm{\sigma_{\mathcal{T}}}.$

\medskip\noindent
The bounds of $\Xi_{i},~i=1,2,3$ in  \eqref{estm:sigma} with $s>\frac{1}{2}$ yield
\begin{align}\label{estm:007_19}
 \norm{\sigma_{\mathcal{T}}}&\lesssim h^{s} ( \sum_{T\in \mathcal{T}_h} (\norm{\nabla(y-\mathcal{E}_{\mathcal{T}}^{k+1}y)}_T^2 + h_T\norm{\nabla(y-\mathcal{E}_{\mathcal{T}}^{k+1}y)}_{\partial T}^2))^{1/2}\nonumber\\&+h_{T}^{1+s} \lVert (f+u)-\Pi_{\mathcal{T}}^{k}(f+u)\rVert,   
\end{align}
 Use the regularity $y\in H^{k+2}(\mathcal{T}_h)$ and $f\in H^{k}(\mathcal{T}_h)$ from the Lemma \ref{lem:regularity_UC} the above inequality to obtain 
\begin{align*}
\norm{\sigma_{\mathcal{T}}}=\norm{\Pi^k_{\mathcal{T}}y-\tilde{y}_{\mathcal{T}}(u)}\lesssim h^{k+1+s}.  
\end{align*}
A combination of steps 1-2 concludes the lemma.
\end{proof}
Similar arguments can be used to prove the following lemma.
\begin{lemma}\label{lem:estm_aux_adj_state}
Under the regularity of solution in Lemma \ref{lem:regularity_UC}, we have the following convergence
   $$\norm{\phi-\tilde{\phi}_{\mathcal{T}}(u)}= \mathcal{O}(h^{k+1}).$$ 
\end{lemma}
\begin{lemma}[$L^2$-error estimates of control]\label{thm:control_estm}
   Under the regularity of solution in Lemma \ref{lem:regularity_UC}, we have the following convergence result for the control variable
    \begin{align*}
     \norm{u-u_h}=\mathcal{O}(h^{k+1}).   
    \end{align*}
\end{lemma}
\begin{proof}
  We recall Theorem \ref{thm:L2_control_error} for the estimate of $\norm{u-u_h}$.
  \begin{align*}
     \norm{u-u_h} \leq &C(C_p,\lambda,\alpha_0) \big(\norm{u - v_h} + \norm{y - \tilde{y}_{\mathcal{T}} (u)}\\&+ \norm{\phi -\tilde{\phi}_{\mathcal{T}} (u)}\big)\quad \forall v_h \in U_h,   
\end{align*}
  Therefore it is enough to estimate $\norm{u - \Pi_h^ku}, \norm{y - \tilde{y}_{\mathcal{T}} (u)},$ and $\norm{\phi -\tilde{\phi}_{\mathcal{T}} (u)}$, where $\Pi_h^ku\in U_h$ be the $L^2$- orthogonal projection of $u.$ Since $u\in H^{k+2}(\mathcal{T}_h)$ we have the following estimate $\norm{u - \Pi^k_hu}\leq h^{k+1} \norm{u}_{H^{k+2}(\mathcal{T}_h)}.$ We have the estimates of $\norm{y - \tilde{y}_{\mathcal{T}} (u)}$ and $\norm{\phi - \tilde{\phi}_{\mathcal{T}} (u)}$ from Lemmas \ref{lem:estm_aux_state} and \ref{lem:estm_aux_adj_state}. Thus we obtain the desired estimate.

\end{proof}

\begin{theorem}[energy error estimate of state]\label{thm:energy_estm_state}
There holds 
\begin{align*}
    \norm{\widehat{I}_h y - \widehat{y}_h}_{1,h} \lesssim ( \sum_{T\in \mathcal{T}_h} (\norm{\nabla(y-\mathcal{E}_{\mathcal{T}}^{k+1}y)}_T^2 + h_T\norm{\nabla(y-\mathcal{E}_{\mathcal{T}}^{k+1}y)}_{\partial T}^2))^{1/2} + \norm{u-u_h}.
\end{align*}
Furthermore, the regularity of solution from Lemma \ref{lem:regularity_UC} leads to
\begin{align*}
	\norm{\widehat{I}_h y - \widehat{y}_h}_{1,h}=\mathcal{O}({h^{k+1}}).	
	\end{align*}
\end{theorem}
\begin{proof}
Utilize triangle inequality to get
\begin{align}\label{auxestm:08}
\norm{\widehat{I}_h y - \widehat{y}_h}_{1,h}&\leq  \norm{\widehat{I}_h y - \tilde{y}_h(u)}_{1,h}+\norm{\tilde{y}_h(u) - \widehat{y}_h}_{1,h}.
\end{align}
{\it Control for $\norm{\tilde{y}_h(u)-\widehat{y}_h}_{1,h}$}: 
A subtraction of \eqref{eqn:disc_state} from \eqref{eqn:aux_y(u)_UM1} yields 
\begin{align}\label{eqn:aux06}
        a_h(\widehat{y}_h-\tilde{y}_h(u), \widehat{w}_h) = (u_h-u, w_{\mathcal{T}}) \quad \forall \;\widehat{w}_h \in \widehat{V}^k_{h,0}.
\end{align}
The choice $\widehat{w}_h=\widehat{y}_h-\tilde{y}_h(u)$ in above equation
and Poincar\'e inequality establish
\begin{align}\label{auxestm:07}
\norm{\tilde{y}_h(u) - \widehat{y}_h}_{1,h}\leq \alpha_0^{-1} C_p \norm{u-u_h}.  
 \end{align}
 We estimate $\norm{\tilde{y}_h(u)-\widehat{I}_h y}_{1,h}$ as follows: Denote $\hhoelem{\sigma} := \tilde{y}_h(u) - \ih y,$ The ellipticity of the bilinear form $a_h(\bullet,\bullet)$ (see \cite[Lemma 2.6]{Book:HHO_Ern:2021}), \eqref{eqn:aux_y(u)_UM1}, \eqref{eqn:conti_state}, and \eqref{eqn:conti_state} reveal
\begin{align}
    \alpha_0 \hhonorm{\hhoelem{\sigma}}^2 &\leq a_h(\hhoelem{\sigma}, \hhoelem{\sigma}) = a_h(\tilde{y}_h(u)-\ih y, \hhoelem{\sigma})= a_h(\tilde{y}_h(u), \hhoelem{\sigma}) - a_h(\ih y, \hhoelem{\sigma})\nonumber\\
    &\leq (f+u,\sigma_{\mathcal{T}}) - a_h(\ih y, \hhoelem{\sigma})
    \leq (-\Delta y,\sigma_{\mathcal{T}}) - a_h(\ih y, \hhoelem{\sigma}).
    \end{align}
The integration by parts formula \eqref{int_part_2}, the definitions of $a_{h}(\bullet,\bullet)$ and $\mathcal{R}_{\mathcal{T}} \circ \widehat{I}_{h}^{k}=\mathcal{E}_{\mathcal{T}}^{k+1}$, and elementary manipulations reveal
   \begin{align}
    \alpha_0 \hhonorm{\hhoelem{\sigma}}^2  &\leq \sum_{T\in \mathcal{T}_h} ((\nabla y, \nabla \sigma_T)_{T}-(\nabla y \cdot \mathbf{n}_T, \sigma_T)_{L^2(\partial T)}) - a_h(\ih y, \hhoelem{\sigma}) \nonumber\\
    &\leq \sum_{T\in \mathcal{T}_h} ((\nabla y, \nabla \sigma_T)_{T}-(\nabla y \cdot \mathbf{n}_T, \sigma_T-\sigma_{\partial T})_{L^2(\partial T)} - (\nabla \mathcal{E}^{k+1}_T y, \nabla \sigma_T)_{T}.\nonumber\\&+(\nabla \mathcal{E}^{k+1}_T y \cdot \mathbf{n}_T, \sigma_T -\sigma_{\partial T})_{L^2(\partial T)} -S_T(\ih y, \hhoelem{\sigma})) \nonumber\\
     &\leq -\sum_{T\in \mathcal{T}_h} ((\nabla (y-\mathcal{E}^{k+1}_T y) \cdot \mathbf{n}_T, \sigma_T-\sigma_{\partial T})_{L^2(\partial T)} +S_T(\ih y, \hhoelem{\sigma})).\nonumber
\end{align}
The Cauchy-Schwarz inequality, Lemma 2.7 \& 2.9 from \cite{Book:HHO_Ern:2021}, establish
\begin{align}\label{auxestm:09}
\norm{\widehat{I}_h y - \tilde{y}_h(u)}_{1,h} \lesssim ( \sum_{T\in \mathcal{T}_h} (\norm{\nabla(y-\mathcal{E}_{\mathcal{T}}^{k+1}y)}_T^2 + h_T\norm{\nabla(y-\mathcal{E}_{\mathcal{T}}^{k+1}y)}_{\partial T}^2))^{1/2}.  
\end{align}
Utilize \eqref{auxestm:07} and \eqref{auxestm:09} in \eqref{auxestm:08} to obtain the first estimate. For the second estimate, we have the following convergence result of elliptic projection 
 $$( \sum_{T\in \mathcal{T}_h} (\norm{\nabla(y-\mathcal{E}_{\mathcal{T}}^{k+1}y)}_T^2 + h_T\norm{\nabla(y-\mathcal{E}_{\mathcal{T}}^{k+1}y)}_{\partial T}^2))^{1/2}\lesssim h^{k+1},$$
 and the estimate of $\norm{u-u_h}$ is available in Lemma \ref{thm:control_estm}. Thus we prove the last estimate. 
\end{proof}
\begin{theorem}[energy error estimate of adjoint state]\label{thm:energy_estm_adj_state_S1UC}
There holds 
\begin{align*}
    \norm{\widehat{I}_h \phi - \widehat{\phi}_h}_{1,h} \lesssim & ( \sum_{T\in \mathcal{T}_h} (\norm{\nabla(\phi-\mathcal{E}_{\mathcal{T}}^{k+1}\phi)}_T^2 + h_T\norm{\nabla(\phi-\mathcal{E}_{\mathcal{T}}^{k+1}\phi)}_{\partial T}^2))^{1/2} \\&+ \norm{u-u_h}+\norm{y-\tilde{y}_{\mathcal{T}}(u)}.
\end{align*} 
Furthermore,  the regularity of solution from Lemma \ref{lem:regularity_UC} yields
\begin{align*}
	\norm{\widehat{I}_h \phi - \widehat{\phi}_h}_{1,h}\lesssim h^{k+1}.	
	\end{align*}
\end{theorem}
\begin{proof}
Introduce the split to get
    \begin{align}\label{auxestm:08_1}
       \norm{\widehat{I}_h \phi - \widehat{\phi}_h}_{1,h}&\leq  \norm{\widehat{I}_h \phi - \tilde{\phi}_h(u)}_{1,h}+\norm{\tilde{\phi}_h(u) - \widehat{\phi}_h}_{1,h}.
    \end{align}
{\it Control for $\norm{\tilde{\phi}_h(u) - \widehat{\phi}_h}_{1,h}$}: Subtract \eqref{eqn:disc_adjstate} from \eqref{eqn:aux_p(u)_UM1} to obtain
\begin{align*}
a_h(\widehat{w}_h, \widehat{\phi}_h-\tilde{\phi}_h(u)) = (y_{\mathcal{T}}-y, w_{\mathcal{T}}) \quad \forall \widehat{w}_h \in \widehat{V}^k_{h,0}.
\end{align*}
The choice $\widehat{w}_h = \widehat{\phi}_h-\tilde{\phi}_h(u)$ in above equation
Poincar\'e inequality and Cauchy-Schwarz inequality yield
\begin{align}\label{auxestm:07_1}
          \alpha_0\norm{\tilde{\phi}_h(u) - \widehat{\phi}_h}_{1,h}&\leq C_p \norm{y-y_{\mathcal{T}}}.  
\end{align}
We estimate $\norm{y-y_{\mathcal{T}}}$ as follows: First introduce the split to get
    \begin{align}\label{estm_007}
        \norm{y-y_{\mathcal{T}}} \leq \norm{y-\tilde{y}_{\mathcal{T}}(u)} + \norm{\tilde{y}_{\mathcal{T}}(u)-y_{\mathcal{T}}}. 
    \end{align}
    Use \eqref{auxestm:07}, with Poincar\'e inequality in the equation \eqref{estm_007} to obtain
    \begin{align}\label{estm:007_1}
        \norm{y-y_{\mathcal{T}}} \leq C(C_p,\alpha_0) (\norm{y-\tilde{y}_{\mathcal{T}}(u)} + \norm{u-u_h}). 
    \end{align}
The combination of \eqref{auxestm:07_1} and \eqref{estm:007_1} yields
 \begin{align}\label{estm:007_11}
   \norm{\tilde{\phi}_h(u) - \widehat{\phi}_h}_{1,h}  \leq C(C_p,\alpha_0) (\norm{y-\tilde{y}_{\mathcal{T}}(u)} + \norm{u-u_h}).  
 \end{align}
 The estimate of $\norm{\widehat{I}_h \phi - \tilde{\phi}_h(u)}_{1,h}$ follows similar to the estimate of $\norm{\widehat{I}_h y - \tilde{y}_h(u)}_{1,h},$ (see the proof of Theorem \ref{thm:energy_estm_state}), and it is of the form:
 \begin{align}\label{auxestm:09_1}
   \norm{\widehat{I}_h \phi - \tilde{\phi}_h(u)}_{1,h} \lesssim ( \sum_{T\in \mathcal{T}_h}( \norm{\nabla(\phi-\mathcal{E}_{\mathcal{T}}^{k+1}\phi)}_T^2 + h_T\norm{\nabla(\phi-\mathcal{E}_{\mathcal{T}}^{k+1}\phi)}_{\partial T}^2))^{1/2}.  
 \end{align}
 Use \eqref{estm:007_11}, \eqref{auxestm:07_1}, and \eqref{auxestm:09_1} in \eqref{auxestm:08_1} we obtain the first estimate. For the second estimate, we have the following convergence result of elliptic projection 
 $$( \sum_{T\in \mathcal{T}_h} (\norm{\nabla(\phi-\mathcal{E}_{\mathcal{T}}^{k+1}\phi)}_T^2 + h_T\norm{\nabla(\phi-\mathcal{E}_{\mathcal{T}}^{k+1}\phi)}_{\partial T}^2))^{1/2}\lesssim h^{k+1},$$
 and the estimate of $\norm{u-u_h}$ and $\norm{y-\tilde{y}_{\mathcal{T}}(u)}$ are available in Lemma \ref{thm:control_estm} and Lemma \ref{lem:estm_aux_state}. Thus we prove the last estimate. 
\end{proof}

\subsection{Scheme UC2 (Variational Discretization):}\label{Scheme:UC_2}
In this scheme, we consider the case where we do not discretize the control variable; we only discretize the state and adjoint state variables. We propose the discrete control problem as follows:
\begin{equation}\label{eqn:disc_costfun_UCVD}
	\min J(\widehat{y}_h, v)=\frac{1}{2}\lVert y_{\mathcal{T}}-y_{d}\rVert^{2}+\frac{\lambda}{2}\|u\|^{2}, 
\end{equation} 
subject to the PDE,
\begin{equation}\label{eqn:disc_state0_UCVD}
	a_h(\widehat{y}_h, \widehat{w}_h)=(f,w_{\mathcal{T}})+(u, w_{\mathcal{T}}) \quad \widehat{w}_h \in \widehat{V}_{h,0}^{k}.    
\end{equation}
Controls $u \in L^2(\Omega).$ The proof of the existence and uniqueness of the problem \eqref{eqn:disc_costfun_UCVD}-\eqref{eqn:disc_state0_UCVD} follows similarly to the proof of Theorem \ref{thm:existance_thm}. 

\noindent
The discrete optimality system of the problem \eqref{eqn:disc_costfun_UCVD}-\eqref{eqn:disc_state0_UCVD}, reads as: find $(\widehat{y}_h, u_h, \widehat{\phi}_h)\in \widehat{V}_{h,0}^{k} \times U_{ad} \times \widehat{V}_{h,0}^{k} $ s.t.
\begin{align}
	a_h(\widehat{y}_h, \widehat{w}_h)&=(f, w_{\mathcal{T}})+(u_h, w_{\mathcal{T}}) \quad \widehat{w}_h \in \widehat{V}_{h,0}^k,\label{eqn:disc_state_UCVD}\\
	a_h(\widehat{w}_h, \widehat{\phi}_h)&=(y_{\mathcal{T}}, w_{\mathcal{T}})-(y_d, w_{\mathcal{T}} ) \quad \forall \widehat{w}_h \in \widehat{V}_{h, 0}^k,\label{eqn:disc_adjstate_UCVD}\\
u_h &= -\frac{\phi_{\mathcal{T}}}{\lambda}.  \label{eqn:disc_VI_UCVD} 
\end{align}

\begin{lemma}[$L^2$-error estimate of control]\label{thm:control_conv_rates_UCVD}
	Let $u$ be the continuous optimal control and $u_h$ be the discrete optimal control solves \eqref{eqn:disc_state_UCVD}-\eqref{eqn:disc_VI_UCVD}. Then the following holds
	\begin{align*}
		\norm{u-u_h}  &\lesssim \norm{y - \tilde{y}_{\mathcal{T}}(u)}+ \norm{\phi - \tilde{\phi}_{\mathcal{T}}(u)},
	\end{align*}
	where $\tilde{\phi}_h(u), \tilde{y}_h  (u)$ solves \eqref{eqn:aux_y(u)_UM1} and \eqref{eqn:aux_p(u)_UM1} respectively. Furthermore, the regularity of solution from Lemma \ref{lem:regularity_UC} leads to
 \begin{align*}
	\norm{u-u_h} =\mathcal{O} (h^{k+1}).
\end{align*}
\end{lemma}

\begin{proof} The notion of $j_{h}$ follows from scheme UC1. An application of second-order sufficient condition \eqref{eqn:disc_red_costfun_UM1} yield
\begin{align*}
\norm{u-u_h}^2&\leq j_h''(w)(u-u_h)^2 = j_h'(u)(u-u_h) - j_h'(u_h)(u-u_h)\nonumber\\&
 \leq (\phi_{\mathcal{T}}(u)+\lambda u, u-u_h) \leq (\phi_{\mathcal{T}}(u)-\phi, u-u_h).\nonumber
 \end{align*}
 Use Cauchy-Schwarz inequality on the right-hand side of the above result to establish
 \begin{align*}
\norm{u-u_h} &\leq \norm{\phi - \phi_{\mathcal{T}}(u)}.\label{estm:aux50}
\end{align*}
Introduce the split to get
$\norm{\phi - \phi_{\mathcal{T}}(u)} \leq \norm{\phi - \tilde{\phi}_{\mathcal{T}}(u)}+\norm{\tilde{\phi}_{\mathcal{T}}(u)- \phi_{\mathcal{T}}(u)} \lesssim \norm{\phi - \tilde{\phi}_{\mathcal{T}}(u)}+\norm{y-\tilde{y}_{\mathcal{T}}(u)}$
where the last inequality follows from the stability estimate of the following auxiliary equation
$a_h(\hhoelem{w},\tilde{\phi}_h(u)-\hhoelem{\phi}(u))= (y-\tilde{y}_{\mathcal{T}}(u),w_{\mathcal{T}})\quad \forall\;\hhoelem{w}\in \Vsp.$
All this, with the above-displayed result, leads to
\begin{align*}
	\norm{u-u_h} \lesssim \norm{\phi - \tilde{\phi}_{\mathcal{T}}(u)}+\norm{y-\tilde{y}_{\mathcal{T}}(u)}.
\end{align*}
This completes the proof of the first part of the lemma. The second part follows from Lemmas \ref{lem:estm_aux_state} -\ref{lem:estm_aux_adj_state}.
\end{proof}
\noindent Proofs of the following two theorems on the energy error estimates for the state and adjoint state are analogous to Theorems 
\ref{lem:estm_aux_state}- \ref{lem:estm_aux_adj_state}.
\begin{theorem}[energy error estimate of state]
There holds 
\begin{align*}
    \norm{\widehat{I}_h y - \widehat{y}_h}_{1,h} \lesssim ( \sum_{T\in \mathcal{T}_h} (\norm{\nabla(y-\mathcal{E}_{\mathcal{T}}^{k+1}y)}_T^2 + h_T\norm{\nabla(y-\mathcal{E}_{\mathcal{T}}^{k+1}y)}_{\partial T}^2))^{1/2} + \norm{u-u_h}.
\end{align*}
Furthermore, the regularity of solution from Lemma \ref{lem:regularity_UC} leads to
\begin{align*}
	\norm{\widehat{I}_h y - \widehat{y}_h}_{1,h} = \mathcal{O}(h^{k+1}).	
	\end{align*}
\end{theorem}
\begin{theorem}[energy error estimate of adjoint state]\label{thm:energy_estm_adj_state}
There holds 
\begin{align*}
    \norm{\widehat{I}_h \phi - \widehat{\phi}_h}_{1,h} \lesssim & ( \sum_{T\in \mathcal{T}_h} (\norm{\nabla(\phi-\mathcal{E}_{\mathcal{T}}^{k+1}\phi)}_T^2 + h_T\norm{\nabla(\phi-\mathcal{E}_{\mathcal{T}}^{k+1}\phi)}_{\partial T}^2))^{1/2} \\&+ \norm{u-u_h}+\norm{y-\tilde{y}_{\mathcal{T}}(u)}.
\end{align*} 
Furthermore, the regularity of solution from Lemma \ref{lem:regularity_UC} leads to
\begin{align*}
	\norm{\widehat{I}_h \phi - \widehat{\phi}_h}_{1,h} =\mathcal{O} (h^{k+1}).	
	\end{align*}
\end{theorem}

\subsection{Scheme UC3.1 (Full reconstruction)}\label{Scheme:UC_3_1} In this scheme, we consider the following discrete counterpart of the continuous problem \eqref{eqn:conti_costfun}-\eqref{eqn:conti_state} as
 \begin{equation}\label{eqn:disc_costfun_UM2a}
 \min J(\widehat{y}_h, \mathcal{R}_{\mathcal{T}} \widehat{u}_h)=\frac{1}{2}\|\mathcal{R}_{\mathcal{T}} \widehat{y}_h-y_{d}\|^{2}+\frac{\lambda}{2}\|\mathcal{R}_{\mathcal{T}} \widehat{u}_h\|^{2}, 
\end{equation} 
subject to the PDE,
\begin{equation}\label{eqn:disc_state_UM2a}
 a_h(\widehat{y}_h, \widehat{w}_h)=(f,\mathcal{R}_{\mathcal{T}}\widehat{w}_h)+(\mathcal{R}_{\mathcal{T}}\widehat{u}_h, \mathcal{R}_{\mathcal{T}}\widehat{w}_h) \quad \forall\; \widehat{w}_h \in \widehat{V}_{h,0}^{k}.    
\end{equation}
The test and trial space in this scheme is $\widehat{V}_{h,0}^{k},$ for $k=0,1$ only. The control space $U_h = \{\mathcal{R}_{\mathcal{T}} \widehat{u}_h \mid \widehat{u}_h \in \widehat{V}_{h}^{k}\; \text{with}\; k=0,1\}.$
Let $\widehat{\mathfrak{S}}_h:L^2(\Omega)\times U_h \rightarrow\widehat{V}^k_{h,0}$ be the solution operator of the PDE \eqref{eqn:disc_state_UM2a} defined by $\widehat{\mathfrak{S}}_h(f,\mathcal{R}_{\mathcal{T}} \widehat{u}_h) = \widehat{y}_h$. The reduced cost functional $j_h$ is defined as follows:
    \begin{align*}
        j_h(\ru{u}) & := J_h(\widehat{\mathfrak{S}}_h(f, \ru{u}),\ru{u})\\
        & = \frac{1}{2} \norm{\mathcal{R}_{\mathcal{T}} \widehat{\mathfrak{S}}_h(f,\ru{u})-y_d}^2 +\frac{\lambda}{2} \norm{\ru{u}}^2.
    \end{align*}
Thus we can reformulate the discrete problem \eqref{eqn:disc_costfun_UM2a}-\eqref{eqn:disc_state_UM2a} as follows:
\begin{align}\label{eqn:reduced_prob_UM2a}
  \min_{\ru{u} \in U_h} j_h(\ru{u})=\frac{1}{2} \norm{\mathcal{R}_{\mathcal{T}} \widehat{\mathfrak{S}}_h(f,\ru{u})-y_d}^2 +\frac{\lambda}{2} \norm{\ru{u}}^2. 
\end{align}
Now we state and prove the existence and uniqueness of the solution of the problem \eqref{eqn:reduced_prob_UM2a}.
\begin{theorem}\label{thm:exist_unique_UM2a}
    There exists a unique solution of the problem \eqref{eqn:reduced_prob_UM2a}.
\end{theorem}

\begin{proof}
    Clearly, $j_h(\ru{u})\geq 0$, therefore infimum exists, call it say $\alpha = \inf j_h(\ru{u}).$ Let $\ru{u}^n$ be an infimizing sequence such that $j_h(\ru{u}^n)\rightarrow\alpha$ as $n\rightarrow\infty.$ This implies the sequence $\norm{\ru{u}^n}$ is bounded in $L^2(\Omega).$ Then there exists a subsequence still index by $n$ such that $\ru{u}^n\rightarrow\ru{u}^*$ weakly in $U_h$. Note that $\ru{u}^*\in U_h$ since $U_h$ is weakly closed. The reduced cost functional $j_h$ is continuous and convex, hence it is weakly lower semi-continuous. Thus, 
    $$j_h(\ru{u}^*) \leq \underline{\lim} j_h(\ru{u}^n) = \alpha. $$
    Therefore, $\ru{u}^*$ is a minimizer. The uniqueness follows from the strict convexity of the reduced cost functional $j_h.$
    \end{proof}
 \noindent   
The discrete optimality system of \eqref{eqn:disc_costfun_UM2a}-\eqref{eqn:disc_state_UM2a} reads as: find $(\widehat{y}_h, \ru{u}, \widehat{\phi}_h)\in \widehat{V}_{h,0}^{k} \times U_h \times \widehat{V}_{h,0}^{k} $ s.t.
\begin{align}
     a_h(\widehat{y}_h, \widehat{w}_h)&=(f, \ru{w})+(\ru{u}, \ru{w}) \quad \widehat{w}_h \in \widehat{V}_{h,0}^k\label{eqn:disc_state_UM2a2}\\
   a_h(\widehat{w}_h, \widehat{\phi}_h)&=(\ru{y}-y_d, \ru{w}) \quad \forall \widehat{w}_h \in \widehat{V}_{h, 0}^k\label{eqn:disc_adjstate_UM2a}\\
   \ru{u} &=-\frac{1}{\lambda} \ru{\phi}. \label{eqn:disc_VI_UM2a} 
\end{align}
We need following auxiliary solutions for error analysis: Let $\tilde{y}_h(u) \; \text{and}\;  \tilde{\phi}_h(u) \in \widehat{V}_{h, 0}^{k}$ solve the following:
\begin{align}
 a_h(\tilde{y}_h(u), \widehat{w}_h)&=(f, \ru{w} )+(u, \ru{w} ) \quad \forall \widehat{w}_h \in \widehat{V}_{h, 0}^{k},\label{eqn:aux_y(u)_UM2a} \\
a_h(\widehat{w}_h, \tilde{\phi}_h(u))&=(y-y_d, \ru{w} ) \quad \forall \widehat{w}_h \in \widehat{V}_{h, 0}^{k}\label{eqn:aux_p(u)_UM2a}.
\end{align}

\begin{theorem}[$L^2$-abstract estimates of control]
    Let $u$ be the solution of \eqref{eqn:conti_costfun}-\eqref{eqn:conti_state}, and $\ru{u}$ be the solution of \eqref{eqn:reduced_prob_UM2a}. Then
    \begin{align*}
        \norm{u-\ru{u}} \lesssim \norm{u-\mathcal{E}_{\mathcal{T}}^{k+1}u}+\norm{y-\mathcal{R}_{\mathcal{T}}\tilde{y}_h(u)}+\norm{\phi-\mathcal{R}_{\mathcal{T}}\tilde{\phi}_h(u)}.
    \end{align*}
\end{theorem}
\begin{proof}
    The triangle inequality yields,
\begin{align}\label{eqn:aux19}
              \norm{u-\ru{u}} \leq \norm{u-\mathcal{E}_{\mathcal{T}}^{k+1}u}+\norm{\mathcal{E}_{\mathcal{T}}^{k+1}u- \ru{u}}. 
\end{align}
From the definition $\mathcal{E}_{\mathcal{T}}^{k+1} = \mathcal{R}_{\mathcal{T}}\circ \widehat{I}_h,$ we can rewrite the second term on right-hand side as $\norm{\mathcal{E}_{\mathcal{T}}^{k+1}u- \ru{u}} = \norm{\mathcal{R}_{\mathcal{T}}\widehat{I}_h u-\ru{u}}.$ Now, we estimate this term next.
The seond-order optimal conditions \eqref{eqn:disc_red_costfun_UM1} leads to
\begin{align*}
\norm{\mathcal{R}_{\mathcal{T}}\widehat{I}_h u-\ru{u}}^2 &\leq j_h''(\ru{w})(\mathcal{R}_{\mathcal{T}} \widehat{I}_h u -\ru{u})^2 = j_h'(\mathcal{R}_{\mathcal{T}}\widehat{I}_h u)(\mathcal{R}_{\mathcal{T}} \widehat{I}_h u -\ru{u}) - j_h'(\ru{u})(\mathcal{R}_{\mathcal{T}} \widehat{I}_h u -\ru{u})\nonumber\\
    & \leq (\rt \widehat{\phi}_h(\widehat{I}_h u) +\lambda \rt \widehat{I}_h u, \rt \ih u -\ru{u})\nonumber\\&\leq (\rt \widehat{\phi}_h(\widehat{I}_h u) - \phi, \rt(\ih u - \hhoelem{u}))+ \lambda (\rt \ih u -u, \rt(\ih u - \hhoelem{u})). 
\end{align*}
 In the last step, we use  $(\phi + \lambda u, \rt(\ih u - \hhoelem{u}))=0,$ given that $u = -\phi/\lambda.$ An application of  Cauchy-Schwarz inequality in the above result leads to 
\begin{align}\label{eqn:aux16}
\norm{\mathcal{R}_{\mathcal{T}}\widehat{I}_h u-\ru{u}} \leq\norm{\phi - \rt \widehat{\phi}_h(\widehat{I}_h u)}+ \lambda \norm{u - \mathcal{E}_{\mathcal{T}}^{k+1}u}. 
 \end{align}
Introduce the split to get
\begin{align}\label{eqn:aux17}
\norm{\phi - \rt \widehat{\phi}_h(\widehat{I}_h u)} \leq \norm{\phi - \rt \tilde{\phi}_h(u)} + \norm{\rt \tilde{\phi}_h(u) - \rt \widehat{\phi}_h(\widehat{I}_h u)}, 
 \end{align}
 where $\widehat{\phi}_h(\widehat{I}_h u)$ satisfies 
$a_h(\hhoelem{w},\widehat{\phi}_h(\ih u))= (\rt \hhoelem{y}(\ih u)-y_d,\ru{w})\; \forall \; \hhoelem{w}\in \Vsp.$ 
Subtracting \eqref{eqn:aux_p(u)_UM2a} from this to get 
\begin{align*}
    a_h(\hhoelem{w},\tilde{\phi}_h(u)- \widehat{\phi}_h(\ih u)) = (y- \rt \hhoelem{y}(\ih u), \ru{w}) \quad \forall \; \hhoelem{w}\in \Vsp.
\end{align*}
The choice $\hhoelem{w} = \tilde{\phi}_h(u)- \widehat{\phi}_h(\ih u)$ in 
above equation, Cauchy-Schwarz inequality, and Poincar\'e inequality reveal
\begin{align}\label{eqn:aux0014_1}
\norm{\rt \tilde{\phi}_h(u) - \rt \widehat{\phi}_h(\widehat{I}_h u)} \lesssim \norm{y- \rt \hhoelem{y}(\ih u)}.
\end{align}
Triangle inequality in the right-hand side of the above inequality leads to
\begin{align}\label{eqn:aux0014_2}
\norm{y- \rt \hhoelem{y}(\ih u)} &\leq \norm{y- \rt \tilde{y}_h(u)} + \norm{\rt \tilde{y}_h(u)- \rt \hhoelem{y}(\ih u)},
 \end{align}
 where $\hhoelem{y}(\ih u)$ solves 
  $   a_h(\hhoelem{y}(\ih u), \ru{w}) = (f, \rt \hhoelem{w})+ (\mathcal{E}_{\mathcal{T}}^{k+1}u, \ru{w})\quad \forall \; \hhoelem{w}\in \Vsp.
  $
 Utilize this
 and \eqref{eqn:aux_y(u)_UM2a} to estimate $\norm{\rt \tilde{y}_h(u)- \rt \hhoelem{y}(\ih u)}\lesssim \norm{u-\mathcal{E}_{\mathcal{T}}^{k+1}u}.$ Thus we have
 \begin{align}\label{eqn:aux0014}
    \norm{y- \rt \hhoelem{y}(\ih u)} \leq  \norm{y- \rt \tilde{y}_h(u)}  + \norm{u-\mathcal{E}_{\mathcal{T}}^{k+1}u}.
 \end{align}
Combine \eqref{eqn:aux17}-
\eqref{eqn:aux0014} to get
  $   \norm{\phi - \rt \widehat{\phi}_h(\widehat{I}_h u)} \leq \norm{y- \rt \tilde{y}_h(u)}  + \norm{\phi - \rt \tilde{\phi}_h(u)} + \norm{u-\mathcal{E}_{\mathcal{T}}^{k+1}u}.
  $
Use above result 
in \eqref{eqn:aux16} and substitute the result in \eqref{eqn:aux19} conclude the theorem.
\end{proof}

\begin{lemma}\label{lem:aux_recon_state_estm}
    There holds
    \begin{align*}
        \norm{y-\rt \tilde{y}_h(u)} &\lesssim \norm{y-\mathcal{E}_{\mathcal{T}}^{k+1}y}+h^{k+s+1} \norm{u+f}_{H^{k}(\mathcal{T}_h)}\\&+ h^s( \sum_{T\in \mathcal{T}_h} (\norm{\nabla(y-\mathcal{E}_{\mathcal{T}}^{k+1}y)}_T^2 + h_T\norm{\nabla(y-\mathcal{E}_{\mathcal{T}}^{k+1}y)}_{\partial T}^2))^{1/2},
    \end{align*}
     furthermore, the regularity of solution from Lemma \ref{lem:regularity_UC} leads to
    \begin{align*}
        \norm{y-\rt \tilde{y}_h(u)} & = \mathcal{O}  (h^{k+s+1})\;\; \text{for}\; k =0,1\; {\text and}\; s\in(1/2,1].
    \end{align*}
\end{lemma}
\begin{proof}
The triangle inequality and the definition of $\mathcal{E}_{\mathcal{T}}^{k+1}={R}_{\mathcal{T}}\circ\widehat{\mathcal{I}}_h$ reveal
    \begin{equation}\label{estm_001_1}
       \norm{y-\rt \tilde{y}_h(u)} \leq \norm{y-\mathcal{E}_{\mathcal{T}}^{k+1}y}+\norm{\mathcal{E}_{\mathcal{T}}^{k+1}y-\rt \tilde{y}_h(u)}
       \leq \norm{y-\mathcal{E}_{\mathcal{T}}^{k+1}y}+\norm{\mathcal{R}_{\mathcal{T}}(\widehat{\mathcal{I}}_h y- \tilde{y}_h(u))}.
    \end{equation}
  Denote $\hhoelem{\sigma} = \tilde{y}_h(u) - \ih y.$ Then we estimate the second term in the right-hand side of \eqref{estm_001_1} as follows. Note that
      $\norm{\rt \hhoelem{\sigma}} \lesssim h \norm{\nabla \rt \hhoelem{\sigma}} + \norm{\sigma_{\mathcal{T}}}.$
  Utilize the continuity of $a_h(\bullet,\bullet)$, to get $\norm{\nabla \rt \hhoelem{\sigma}}^2\leq a_h(\hhoelem{\sigma},\hhoelem{\sigma}) \lesssim \hhonorm{\hhoelem{\sigma}}^2.$ All this yields 
\begin{align}\label{estm_001_2}
      \norm{\rt \hhoelem{\sigma}} \lesssim h \hhonorm{\hhoelem{\sigma}} + \norm{\sigma_{\mathcal{T}}}.
\end{align}
Therefore, it is enough to estimate $\hhonorm{\hhoelem{\sigma}}$ and $\norm{\sigma_{\mathcal{T}}}.$ We first estimate $\hhonorm{\hhoelem{\sigma}},$  Lemma 2.9 from \cite{Book:HHO_Ern:2021}, the definition of $\hhoelem{\sigma},$ \eqref{eqn:aux_y(u)_UM2a}, and elementary manipulation reveal
\begin{align}
    \alpha_0 \hhonorm{\hhoelem{\sigma}}^2 &\leq a_h(\hhoelem{\sigma}, \hhoelem{\sigma}) = a_h(\tilde{y}_h(u)-\ih y, \hhoelem{\sigma}) \leq a_h(\tilde{y}_h(u), \hhoelem{\sigma}) - a_h(\ih y, \hhoelem{\sigma})\nonumber\\
    &\leq (f+u,\ru{\sigma}) - a_h(\ih y, \hhoelem{\sigma})
    \leq (f+u,\ru{\sigma}-\sigma_{\mathcal{T}})+(f+u,\sigma_{\mathcal{T}}) - a_h(\ih y, \hhoelem{\sigma}).\label{estm:aux32}
\end{align}
Now, we estimate the each term on the right-hand side of \eqref{estm:aux32}.

\medskip \noindent
{\it Control for $(f+u,\ru{\sigma}-\sigma_{\mathcal{T}})$}: The estimate of $\Pi_T^0,$ and the estimate $\norm{\mathcal{R}_{\mathcal{T}}\widehat{\sigma}_T-\sigma_T}_{T}\lesssim h_T\norm{\nabla(\mathcal{R}_{\mathcal{T}}\widehat{\sigma}_T-\sigma_T)}_{T}\lesssim h_T \norm{\widehat{\sigma}_h}_{1,h}$ lead to
\begin{align}
    &(f+u,\ru{\sigma}-\sigma_{\mathcal{T}}) = \sum_{T\in \mathcal{T}_h} (f+u,\mathcal{R}_{\mathcal{T}}\widehat{\sigma}_T-\sigma_T)_{T}
    = \sum_{T\in \mathcal{T}_h} ((f+u) - \Pi_T^0(f+u),\mathcal{R}_{\mathcal{T}}\widehat{\sigma}_T-\sigma_T)_{T}\nonumber\\
    &\lesssim \sum_{T\in \mathcal{T}_h} h_T^{k}\norm{f+u}_{H^{k}(T)} \norm{\mathcal{R}_{\mathcal{T}}\widehat{\sigma}_T-\sigma_T}_{T}
     \lesssim \sum_{T\in \mathcal{T}_h} h_T^{k+1}\norm{f+u}_{H^{k}(T)} \norm{\nabla(\mathcal{R}_{\mathcal{T}}\widehat{\sigma}_T-\sigma_T)}_{T}
     \nonumber\\
    &\lesssim h^{k+1}\norm{f+u}_{H^{k}(\mathcal{T}_h)} \norm{\widehat{\sigma}_h}_{1,h}.\label{estm:aux30}
\end{align}
{\it Control for $(f+u,\sigma_{\mathcal{T}}) - a_h(\ih y, \hhoelem{\sigma})$}: Use \eqref{eqn:conti_state} and integration by part formula \eqref{int_part_2} to obtain
$(f+u,\sigma_{\mathcal{T}})-a_h(\ih y, \hhoelem{\sigma})= (-\Delta y,\sigma_{\mathcal{T}}) - a_h(\ih y, \hhoelem{\sigma})= \sum_{T\in \mathcal{T}_h} ((\nabla y, \nabla \sigma_T)_{T}-(\nabla y \cdot \mathbf{n}_T, \sigma_T)_{L^2(\partial T)}) - a_h(\ih y, \hhoelem{\sigma})=
    \sum_{T\in \mathcal{T}_h} ((\nabla y, \nabla \sigma_T)_{T}-(\nabla y \cdot \mathbf{n}_T, \sigma_T-\sigma_{\partial T})_{L^2(\partial T)} - (\nabla \mathcal{E}^{k+1}_T y, \nabla \sigma_T)_{T}+(\nabla \mathcal{E}^{k+1}_T y \cdot \mathbf{n}_T, \sigma_T -\sigma_{\partial T})_{L^2(\partial T)} -S_T(\ih y, \hhoelem{\sigma})),$ where definition of  $a_h(\bullet,\bullet)$ and $\mathcal{E}^{k+1}_T=\mathcal{R}_{T}\circ \widehat{I}_{T}^{k}$ used in the last step. Now, use Cauchy-Schwarz inequality, trace inequality \eqref{trace_ineq}, and the definition of $\hhonorm{\bullet}$ from Section \ref{sec_23} yield
   \begin{align} 
    (f+u,\sigma_{\mathcal{T}}) - a_h(\ih y, \hhoelem{\sigma})&\lesssim ( \sum_{T\in \mathcal{T}_h} (\norm{\nabla(y-\mathcal{E}_{\mathcal{T}}^{k+1}y)}_T^2 + h_T\norm{\nabla(y-\mathcal{E}_{\mathcal{T}}^{k+1}y)}_{\partial T}^2))^{1/2}\hhonorm{\hhoelem{\sigma}}.\label{estm:aux31}
\end{align}
The estimates \eqref{estm:aux30} and \eqref{estm:aux31} in \eqref{estm:aux32} give
\begin{align}\label{estm:008_2}
   \hhonorm{\hhoelem{\sigma}} \lesssim&  ( \sum_{T\in \mathcal{T}_h} (\norm{\nabla(y-\mathcal{E}_{\mathcal{T}}^{k+1}y)}_T^2 + h_T\norm{\nabla(y-\mathcal{E}_{\mathcal{T}}^{k+1}y)}_{\partial T}^2))^{1/2}  + h^{k+1} \norm{f+u}_{H^{k}(\mathcal{T}_h)} . 
\end{align}
{\it Control of $\norm{\sigma}_{\mathcal{T}}$}:
Use the similar arguments as in \eqref{estm:007_19} to achieve
\begin{align}\label{estm:008_1}
 \norm{\sigma_{\mathcal{T}}}
 &\lesssim h^{s} ( \sum_{T\in \mathcal{T}_h} (\norm{\nabla(y-\mathcal{E}_{\mathcal{T}}^{k+1}y)}_T^2 + h_T\norm{\nabla(y-\mathcal{E}_{\mathcal{T}}^{k+1}y)}_{\partial T}^2))^{1/2}+h^{1+s} \lVert (f+u)-\Pi_{\mathcal{T}}^{k}(f+u)\rVert,   
\end{align}
 for $s\in (\frac{1}{2},1].$ Utilize the estimates \eqref{estm:008_2}-\eqref{estm:008_1} in \eqref{estm_001_2} to get
 \begin{align}\label{estm:008_3}
     \norm{\rt \hhoelem{\sigma}}&\lesssim  h^s( \sum_{T\in \mathcal{T}_h} (\norm{\nabla(y-\mathcal{E}_{\mathcal{T}}^{k+1}y)}_T^2 + h_T\norm{\nabla(y-\mathcal{E}_{\mathcal{T}}^{k+1}y)}_{\partial T}^2))^{1/2}+ h^{k+s+1} \norm{u+f}_{H^{k}(\mathcal{T}_h)}.
 \end{align}
 Use \eqref{estm:008_3} in \eqref{estm_001_1} to obtain the first estimate. The second estimate follows from invoking the regularity of solution from Lemma \ref{lem:regularity_UC}. This completes the lemma.
\end{proof}

\medskip\noindent
We state the following lemma for the estimate of $\norm{\phi-\rt \tilde{\phi}_h(u)}.$ The proof can be derived similarly to the proof of Lemma \ref{lem:aux_recon_state_estm}. 
\begin{lemma}\label{lem:aux_recon_adj_state_estm}
    There holds
    \begin{align*}
        \norm{\phi-\rt \tilde{\phi}_h(u)} &\lesssim \norm{\phi-\mathcal{E}_{\mathcal{T}}^{k+1}\phi}+h^{k+s+1} \norm{y-y_d}_{H^{k}(\mathcal{T}_h)}\\&+ h^s( \sum_{T\in \mathcal{T}_h} (\norm{\nabla(\phi-\mathcal{E}_{\mathcal{T}}^{k+1}\phi)}_T^2 + h_T\norm{\nabla(\phi-\mathcal{E}_{\mathcal{T}}^{k+1}\phi)}_{\partial T}^2))^{1/2}.
    \end{align*}
        Furthermore, the regularity of solution from Lemma \ref{lem:regularity_UC} leads to
    \begin{align*}
        \norm{\phi-\rt \tilde{\phi}_h(u)} &= \mathcal{O} (h^{k+s+1}) \;\; \text{for}\; k =0,1\; {\text and}\; s\in(1/2,1].
    \end{align*}
\end{lemma}

\medskip\noindent
Utilize the estimates from Lemma \ref{lem:aux_recon_state_estm}- \ref{lem:aux_recon_adj_state_estm} to derive the following error estimate of control.
\begin{lemma}[$L^2$- error estimate of control]\label{lem:control_estm_007}
    Under the regularity of solution in Lemma \ref{lem:regularity_UC}, we have the following convergence result for the control variable
    \begin{align*}
        \norm{u-\ru{u}} = \mathcal{O} (h^{k+s+1}) \;\; \text{for}\; k =0,1\; {\text and}\; s\in(1/2,1].
    \end{align*}
\end{lemma}

\begin{lemma}[Energy error estimate of state]\label{lem:estm_energy_007_1} There holds,
 \begin{align*}
   \hhonorm{\hhoelem{y}-\ih y} \lesssim &( \sum_{T\in \mathcal{T}_h} (\norm{\nabla(y-\mathcal{E}_{\mathcal{T}}^{k+1}y)}_T^2 + h_T\norm{\nabla(y-\mathcal{E}_{\mathcal{T}}^{k+1}y)}_{\partial T}^2))^{1/2}+ h^{k+1} |f+u|_{H^{k}(\mathcal{T}_h)} + \norm{u- \ru{u}}.
 \end{align*}
 Furthermore, the regularity of solution in Lemma \ref{lem:regularity_UC} leads to
\begin{align*}
\hhonorm{\hhoelem{y}-\ih y} = \mathcal{O} (h^{k+1}) \;\; \text{for}\; k =0,1.
\end{align*}
\end{lemma}

\begin{proof} 
Introduce the split to get
\begin{align}\label{estm:aux34}
\hhonorm{\hhoelem{y}-\ih y} \lesssim  \hhonorm{\tilde{y}_h(u)-\ih y} + \hhonorm{\hhoelem{y}-\tilde{y}_h(u)}. 
\end{align}
Now, we estimate each term on the right-hand side of \eqref{estm:aux34}.

\medskip\noindent
{\it Control of $\hhonorm{\tilde{y}_h(u)-\ih y}$}: From \eqref{estm:008_2}, it holds
\begin{align}
\hhonorm{\tilde{y}_h(u)-\ih y} \lesssim&  ( \sum_{T\in \mathcal{T}_h} (\norm{\nabla(y-\mathcal{E}_{\mathcal{T}}^{k+1}y)}_T^2 + h_T\norm{\nabla(y-\mathcal{E}_{\mathcal{T}}^{k+1}y)}_{\partial T}^2))^{1/2}+ h^{k+1} |f+u|_{H^{k}(\mathcal{T}_h)}.\label{estm:aux33}
\end{align}
{\it Control of $\hhonorm{\hhoelem{y} - \tilde{y}_h(u)}$}:
 subtraction of \eqref{eqn:aux_y(u)_UM2a} from \eqref{eqn:disc_state_UM2a} leads to the following error equation
 \begin{align}\label{eqn:aux42}
     a_h(\tilde{y}_h(u)-\hhoelem{y},\hhoelem{w}) = (u- \ru{u}, \ru{w}) \quad \hhoelem{w} \in \Vsp
 \end{align}
 The choice $\hhoelem{w} = \tilde{y}_h(u)-\hhoelem{y}$ in \eqref{eqn:aux42}, yields $\hhonorm{\tilde{y}_h(u)-\hhoelem{y}}\lesssim \norm{u- \ru{u}}.$ Use this estimate and \eqref{estm:aux33} in \eqref{estm:aux34} to derive 
 \begin{align*}
   \hhonorm{\hhoelem{y}-\ih y} \lesssim \big( \sum_{T\in \mathcal{T}_h} (\norm{\nabla(y-\mathcal{E}_{\mathcal{T}}^{k+1}y)}_T^2 + h_T\norm{\nabla(y-\mathcal{E}_{\mathcal{T}}^{k+1}y)}_{\partial T}^2)\big)^{1/2} + h^{k+1} |f+u|_{H^{k}(\mathcal{T}_h)} + \norm{u- \ru{u}}.
 \end{align*}
 The second estimate follows from the regularity result from Lemma \ref{lem:regularity_UC} and Lemma \ref{lem:control_estm_007}.
\end{proof}

\medskip \noindent
We state the following lemma for the energy error estimate of the adjoint state. The proof can be derived similarly to the proof of Lemma \ref{lem:estm_energy_007_1}.

\begin{theorem}[energy error estimate of adjoint state] There holds,
 \begin{align*}
   \hhonorm{\hhoelem{\phi}-\ih \phi} \lesssim &( \sum_{T\in \mathcal{T}_h} (\norm{\nabla(\phi-\mathcal{E}_{\mathcal{T}}^{k+1}\phi)}_T^2 + h_T\norm{\nabla(\phi-\mathcal{E}_{\mathcal{T}}^{k+1}\phi)}_{\partial T}^2))^{1/2}\\&+ h^{k+1} |y-y_d|_{H^{k}(\mathcal{T}_h)}+ \norm{y- \ru{y}} + \norm{u- \ru{u}}.
 \end{align*}  
 Furthermore, the regularity of solution in Lemma \ref{lem:regularity_UC} leads to
  \begin{align*}
        \hhonorm{\hhoelem{\phi}-\ih \phi}= \mathcal{O}  (h^{k+1}),\;\; \text{for}\;k =0,1.
    \end{align*}
\end{theorem}

\subsection{Scheme UC3.2 (Partial reconstruction)}\label{Scheme:UC_3_2} In this method, we consider the following discrete counterpart of the continuous problem \eqref{eqn:conti_costfun}-\eqref{eqn:conti_state} as
 \begin{equation}\label{eqn:disc_costfun_UM2b}
 \min J(\widehat{y}_h, \mathcal{R}_{\mathcal{T}} \widehat{u}_h)=\frac{1}{2}\|y_{\mathcal{T}}-y_{d}\|^{2}+\frac{\lambda}{2}\|\mathcal{R}_{\mathcal{T}} \widehat{u}_h\|^{2}, 
\end{equation} 
subject to the PDE,
\begin{equation}\label{eqn:disc_state_UM2b1}
 a_h(\widehat{y}_h, \widehat{w}_h)=(f,w_{\mathcal{T}})+(\mathcal{R}_{\mathcal{T}}\widehat{u}_h, w_{\mathcal{T}}) \quad \widehat{w}_h \in \widehat{V}_{h,0}^{k+}.    
\end{equation}
The control space $U_h = \{\mathcal{R}_{\mathcal{T}} \widehat{u}_h \mid \widehat{u}_h \in \widehat{V}_{h}^{k}\}.$ In this scheme, we always assume $k\geq 2.$
 The proof of existence and uniqueness follows similarly to the proof of Theorem \ref{thm:exist_unique_UM2a}. The discrete optimality system of \eqref{eqn:disc_costfun_UM2b}-\eqref{eqn:disc_state_UM2b1} reads as: find $(\widehat{y}_h, \ru{u}, \widehat{\phi}_h)\in \widehat{V}_{h,0}^{k+} \times U_h \times \widehat{V}_{h,0}^{k+} $ such that
\begin{align}
     a_h(\widehat{y}_h, \widehat{w}_h)&=(f,  w_{\mathcal{T}})+(\ru{u}, w_{\mathcal{T}} ) \quad \widehat{w}_h \in \widehat{V}_{h,0}^{k+}\label{eqn:disc_state_UM2b}\\
   a_h(\widehat{w}_h, \widehat{\phi}_h)&=(\ru{y}-y_d, w_{\mathcal{T}}) \quad \forall \widehat{w}_h \in \widehat{V}_{h, 0}^{k+}\label{eqn:disc_adjstate_UM2b}\\
   \ru{u} &=-\frac{1}{\lambda} \phi_{\mathcal{T}}. \label{eqn:disc_VI_UM2b} 
\end{align}
We need the following auxiliary solutions for error analysis: Let $\tilde{y}_h(u), \tilde{\phi}_h(u) \in \widehat{V}_{h, 0}^{k+}$ solve the following:
\begin{align}
 a_h(\tilde{y}_h(u), \widehat{w}_h)&=(f, w_{\mathcal{T}} )+(u, w_{\mathcal{T}} ) \quad \forall \widehat{w}_h \in \widehat{V}_{h, 0}^{k+},\label{eqn:aux_y(u)_UM2b} \\
a_h(\widehat{w}_h, \tilde{\phi}_h(u))&=(y-y_d, w_{\mathcal{T}} ) \quad \forall \widehat{w}_h \in \widehat{V}_{h, 0}^{k+}\label{eqn:aux_p(u)_UM2b}.
\end{align}
Below we state the error estimate results for control, state, and adjoint state variables. The proof follows from subsection \ref{Scheme:UC_3_1}.
\begin{theorem}[$L^2$- error estimates of control]
    Let $u$ be the continuous control satisfies \eqref{eqn:conti_os_state}-\eqref{eqn:conti_os_control} and $\ru{u}$ be the discrete control satisfies \eqref{eqn:disc_state_UM2b}-\eqref{eqn:disc_VI_UM2b}. There holds,
    \begin{align*}
        \norm{u-\ru{u}} \lesssim \norm{u-\mathcal{E}_{\mathcal{T}}^{k+1}u}+\norm{y-\tilde{y}_{\mathcal{T}}(u)}+\norm{p-\tilde{p}_{\mathcal{T}}(u)}.
    \end{align*}
 Furthermore, the regularity of solution in Lemma \ref{lem:regularity_UC} leads to
\begin{align*}
\norm{u-\ru{u}} = \mathcal{O} (h^{k+s+1})\;\; \text{for}\; k\geq 2\; \text{and}\;s\in(1/2,1]. 
\end{align*}
\end{theorem}

\begin{theorem}[Energy error estimate of state] The following estimates holds,
 \begin{align*}
   \hhonorm{\hhoelem{y}-\ih y} \lesssim & \big( \sum_{T\in \mathcal{T}_h} (\norm{\nabla(y-\mathcal{E}_{\mathcal{T}}^{k+1}y)}_T^2 + h_T\norm{\nabla(y-\mathcal{E}_{\mathcal{T}}^{k+1}y)}_{\partial T}^2)\big)^{1/2}\\&+ \norm{\nabla_{\mathcal{T}}(y-\Pi_{\mathcal{T}}^{k+1}y)} + \norm{u- \ru{u}}.
 \end{align*}
 Furthermore, the regularity of solution in Lemma \ref{lem:regularity_UC} leads to
\begin{align*}
\hhonorm{\hhoelem{y}-\ih y} = \mathcal{O}  (h^{k+1})\;\; \text{for}\; k\geq 2.
\end{align*}  
\end{theorem}
\begin{theorem}[Energy error estimate of adjoint state] The following estimates holds,
 \begin{align*}
   \hhonorm{\hhoelem{\phi}-\ih \phi} \lesssim &\big( \sum_{T\in \mathcal{T}_h} (\norm{\nabla(\phi-\mathcal{E}_{\mathcal{T}}^{k+1}\phi)}_T^2 + h_T\norm{\nabla(\phi-\mathcal{E}_{\mathcal{T}}^{k+1}\phi)}_{\partial T}^2)\big)^{1/2}\\&+  \norm{\nabla_{\mathcal{T}}(\phi-\Pi_{\mathcal{T}}^{k+1}\phi)} + \norm{y- \tilde{y}_{\mathcal{T}}(u)} + \norm{u- \ru{u}}.
 \end{align*}
 Furthermore, the regularity of solution in Lemma \ref{lem:regularity_UC} leads to
\begin{align*}
\hhonorm{\hhoelem{\phi}-\ih \phi} = \mathcal{O}(h^{k+1})\;\; \text{for}\; k\geq 2.
\end{align*}    
\end{theorem}
\section{Control problem with constraints}\label{sec:OCP_HHO_WC}
In this section, we consider the following distributed optimal control problem with control constraints:
\begin{equation}\label{eqn:conti_costfun_WC}
 \min J(y, u)=\frac{1}{2}\lVert y-y_{d}\rVert^{2}+\frac{\lambda}{2}\|u\|^{2},   
\end{equation}
subject to the PDE,
\begin{equation}\label{eqn:conti_state_WC_srtong_f}
\begin{aligned}
-\Delta y & =f+u \quad \text{in } \Omega, \\
y & =0 \quad \quad \text { on } \partial \Omega.
\end{aligned}
\end{equation}
The control $u$ comes from the set $U_{ad}:= \{ u\in L^2(\Omega)\mid u_a \leq u(x)\leq u_b\; \text{a.e.} \; x \in \Omega \}$.

The existence and uniqueness of the solution of the problem \eqref{eqn:conti_costfun_WC}-\eqref{eqn:conti_state_WC_srtong_f} is standard and can be found in \cite[Theorem 2.14]{trolzstch:2005:Book}. The above problem can be formulated as the following optimality system using first-order necessary optimality condition. Find $(y, u, \phi)\in H^1_0(\Omega) \times  L^2(\Omega)\times H^1_0(\Omega) $ s.t.
    \begin{align}
         a(y, v)&=(f, v)+(u, v) \quad \forall v \in H_{0}^1(\Omega),\label{eqn:conti_state_WC} \\ a(v, \phi)&=(y-y_d, v)  \quad\forall v \in H^1_{0}(\Omega),\label{eqn:conti_adjstate_WC}\\ 
        (\phi+\lambda u,  v-u)&\geq 0\quad \forall v \in U_{ad}.\label{eqn:conti_VI_WC}    
    \end{align}
One can rewrite the inequality \eqref{eqn:conti_VI_WC} in a more sophisticated form as
\begin{align}\label{proj_formula}
  u=\mathcal{P}_{U_{ad}}(-\frac{\phi}{\lambda}),  
\end{align}
where $\mathcal{P}_{U_{ad}}(w):=\min\{u_b,\max\{u_a,w\}\}$.
\begin{assumption}\label{Assump_WC}
    Throughout this section, we assume that data $f,y_d\in H^1(\Omega).$
\end{assumption}

\begin{lemma}[Regularity]\label{lem:regularity_WC}
    Under the Assumption \ref{Assump_WC}, the solution of the problem \eqref{eqn:conti_state_WC}-\eqref{eqn:conti_VI_WC} possess the following regularity $y,\phi \in H^3(\Omega),$ and $u \in W^{1,p}(\Omega),\; 2< p < \infty.$(Note that this smoothness holds for both $2$ and $3$ dimensions.)
\end{lemma}

\begin{proof}
    The source term $f\in H^1(\Omega)$ and $u\in L^2(\Omega),$ thus $f+u\in L^2(\Omega),$ hence from the equation \eqref{eqn:conti_state_WC_srtong_f} we have $y\in H^2(\Omega).$ Given that $y_d\in H^1(\Omega)$, thus $y-y_d\in H^1(\Omega),$ therefore from the equation \eqref{eqn:conti_adjstate_WC}, we have $\phi\in H^3(\Omega).$ Therefore $\phi\in W^{1,p}(\Omega),\; 2< p < \infty.$ Now the projection formula \eqref{proj_formula} yields $u\in  W^{1,p}(\Omega),\; 2< p < \infty.$ Thus $f+u\in H^1(\Omega),$ hence $y\in H^3(\Omega).$ Therefore all together $y,\phi \in H^3(\Omega)$ and $u \in W^{1,p}(\Omega),\; 2< p < \infty,$ and this proves the Lemma.
\end{proof}
Following, we propose three discrete schemes to compute and analyze the optimal control problem \eqref{eqn:conti_costfun_WC}-\eqref{eqn:conti_state_WC_srtong_f}.

\subsection{Scheme WC1:}\label{Scheme:WC_1} We consider the following discrete counterpart of the continuous problem \eqref{eqn:conti_costfun_WC}-\eqref{eqn:conti_state_WC_srtong_f} as
 \begin{equation}\label{eqn:disc_costfun_WCM0}
 \min J(\widehat{y}_h, u_h)=\frac{1}{2}\lVert y_{\mathcal{T}}-y_{d}\rVert^{2}+\frac{\lambda}{2}\|u_h\|^{2}, 
\end{equation} 
subject to the PDE,
\begin{equation}\label{eqn:disc_state_WCM01}
 a_h(\widehat{y}_h, \widehat{w}_h)=(f,w_{\mathcal{T}})+(u_h, w_{\mathcal{T}}) \quad \widehat{w}_h \in \widehat{V}_{h,0}^{0}.    
\end{equation}
Controls $u_h \in U^h_{ad}:= \{u_h \in \mathbb{P}_0(\mathcal{T}_h) \mid  u_a \leq u_h|_T\leq u_b,\; \;\forall T \in \mathcal{T}_h \}.$

The existence and uniqueness of the solution of the problem \eqref{eqn:disc_costfun_WCM0}-\eqref{eqn:disc_state_WCM01} follows from Theorem \ref{thm:existance_thm}. The discrete optimality system reads as: find $(\widehat{y}_h, u_h, \widehat{\phi}_h)\in \widehat{V}_{h,0}^{0} \times U^h_{ad} \times \widehat{V}_{h,0}^{0} $ s.t.
\begin{align}
     a_h(\widehat{y}_h, \widehat{w}_h)&=(f, w_{\mathcal{T}})+(u_h, w_{\mathcal{T}}) \quad \forall \widehat{w}_h \in \widehat{V}_{h,0}^0,\label{eqn:disc_state_WCM0}\\
   a_h(\widehat{w}_h, \widehat{\phi}_h)&=(y_{\mathcal{T}}, v_{\mathcal{T}})-(y_d, w_{\mathcal{T}} ) \quad \forall \widehat{w}_h \in \widehat{V}_{h, 0}^0,\label{eqn:disc_adjstate_WCM0}\\
   (\phi_{\mathcal{T}} + \lambda u_h, v_h-u_h)&\geq 0 \quad v_h\in U^h_{ad}.  \label{eqn:disc_VI_WCM0} 
\end{align}
We need the following auxiliary problems for error analysis: Let $\tilde{y}_h(u), \tilde{\phi}_h(u) \in \widehat{V}_{h, 0}^{0}$ solve the following:
\begin{align}
 a_h(\tilde{y}_h(u), \widehat{w}_h)&=(f, w_{\mathcal{T}} )+(u, w_{\mathcal{T}} ) \quad \forall \widehat{w}_h \in \widehat{V}_{h, 0}^{0},\label{eqn:aux_y(u)_WCM0} \\
a_h(\widehat{w}_h, \tilde{\phi}_h(u))&=(y-y_d, w_{\mathcal{T}} ) \quad \forall \widehat{w}_h \in \widehat{V}_{h, 0}^{0}.\label{eqn:aux_p(u)_WCM0}
\end{align}
\begin{lemma}\label{lem:aux_estm01_0}
    The following estimate holds,
    \begin{align*}
        (u-u_h,\tilde{\phi}_{\mathcal{T}}(u)-\phi_{\mathcal{T}}) = (y-y_{\mathcal{T}}, \tilde{y}_{\mathcal{T}}(u)-y) + \norm{y-y_{\mathcal{T}}}^2.
    \end{align*}
\end{lemma}

\begin{proof}
Subtract \eqref{eqn:disc_state_WCM0} from \eqref{eqn:aux_y(u)_WCM0} to get
    \begin{align}\label{eqn:aux22_0}
        a_h(\tilde{y}_h(u)-\hhoelem{y},\hhoelem{w}) = (u-u_h, w_{\mathcal{T}}),\quad \forall \hhoelem{w} \in \Vsp.
    \end{align}
Subtract \eqref{eqn:disc_adjstate_WCM0} from \eqref{eqn:aux_p(u)_WCM0} to obtain
    \begin{align}\label{eqn:aux23_0}
        a_h(\tilde{\phi}_h(u)-\hhoelem{\phi},\hhoelem{w}) = (y-y_{\mathcal{T}}, w_{\mathcal{T}}),\quad \forall \hhoelem{w} \in \Vsp.
    \end{align}
The choice $\hhoelem{w} = \tilde{y}_h(u)-\hhoelem{y}$ in \eqref{eqn:aux23_0} and $\hhoelem{w} = \tilde{\phi}_h(u)-\hhoelem{\phi}$ in \eqref{eqn:aux22_0} and equating them lead to
\begin{align*}
&(u-u_h, \tilde{\phi}_{\mathcal{T}}(u)-\phi_{\mathcal{T}}) = (y-y_{\mathcal{T}}, \tilde{y}_{\mathcal{T}}(u)-y_{\mathcal{T}})\\
& = (y-y_{\mathcal{T}}, \tilde{y}_{\mathcal{T}}(u)-y) +(y-y_{\mathcal{T}}, y-y_{\mathcal{T}})
 = (y-y_{\mathcal{T}}, \tilde{y}_{\mathcal{T}}(u)-y) + \norm{y-y_{\mathcal{T}}}^2.
    \end{align*}
    This completes the lemma.
\end{proof}
\noindent
\begin{lemma}[$L^2$- error estimate of control and state]\label{thm:control_conv_rates_WC}
   Let $(y,u)$ be the optimal state, and control satisfies \eqref{eqn:conti_os_state}-\eqref{eqn:conti_os_control} with the regularity assumption from Lemma \ref{lem:regularity_WC}. The pair $(\hhoelem{y},u_h)$ is the discrete state and optimal control satisfies \eqref{eqn:disc_state_WCM0}-\eqref{eqn:disc_VI_WCM0}. Then, it holds
\begin{align*}
\norm{u-u_h} + \norm{y-y_h} = \mathcal{O} (h).
\end{align*}
\end{lemma}
\begin{proof}
The choice $v=u_h$ in \eqref{eqn:conti_VI_WC} leads to
        $(\phi+\lambda u,  u_h-u)\geq0.$
    Use \eqref{eqn:disc_VI_WCM0} to get
        $(\phi_{\mathcal{T}} + \lambda u_h, u-u_h) \geq -(\phi_{\mathcal{T}} + \lambda u_h, v_h-u)\quad \forall v_h\in U^h_{ad}.$
    Add these two results and elementary manipulations lead to
$ \lambda \norm{u-u_h}^2 \leq (\phi_{\mathcal{T}} -\phi, u-u_h)+(\phi_{\mathcal{T}} +\lambda u_h, v_h-u)\leq (\phi_{\mathcal{T}} - \tilde{\phi}_{\mathcal{T}}(u), u-u_h)+(\tilde{\phi}_{\mathcal{T}}(u)-\phi, u-u_h)+ (\phi_{\mathcal{T}}-\phi+\lambda(u_h-u),v_h-u)+(\phi+\lambda u, v_h-u),$ where splits are used on the right-hand side in the last step. 
 \begin{align*}  
 \lambda \norm{u-u_h}^2
&\leq (y-\tilde{y}_{\mathcal{T}}(u),y-y_{\mathcal{T}})-\norm{y-y_{\mathcal{T}}}^2
+(\tilde{\phi}(u)-\phi, u-u_h) + (\phi_{\mathcal{T}} -\phi, v_h-u)\\&\qquad +\lambda(u_h-u,v_h-u)+(\phi+\lambda u, v_h-u).
\end{align*}
The elementary manipulations and Cauchy-Schwarz inequality in the above-inequality reveal 
\begin{align*}
\norm{u-u_h}^2+\norm{y-y_{\mathcal{T}}}^2 &\lesssim \norm{y-\tilde{y}_{\mathcal{T}}(u)}^2 + \norm{\phi-\tilde{\phi}(u)}^2 + \norm{u-v_h}^2+(\phi_{\mathcal{T}} -\phi, v_h-u) +(\phi+\lambda u, v_h-u). 
\end{align*}
    Now $\norm{\phi - \phi_{\mathcal{T}}} \lesssim \norm{\phi-\tilde{\phi}_{\mathcal{T}}(u)} + \norm{y-y_{\mathcal{T}}}.$ Therefore, for all $v_h \in U^h_{ad},$
    \begin{align}\label{Eqn:102}
         \norm{u-u_h}^2+\norm{y-y_{\mathcal{T}}}^2& \lesssim \norm{y-\tilde{y}_{\mathcal{T}}(u)}^2 + \norm{\phi-\tilde{\phi}_{\mathcal{T}}(u)}^2 + \norm{u-v_h}^2+(\phi+\lambda u, v_h-u). 
    \end{align}
The choice $v_h = \Pi_h^0 u \in U_h^{ad}$ in \eqref{Eqn:102}, where $\Pi_h^0: L^2(\Omega)\rightarrow\mathbb{P}_0(\mathcal{T}_h)$ be the $L^2$-orthogonal projection onto $\mathbb{P}_0(\mathcal{T}_h).$ Use the regularity of solution from Lemma \ref{lem:regularity_WC} to obtain
    \begin{align*}
    &\norm{u-u_h}^2+\norm{y-y_{\mathcal{T}}}^2 \lesssim \norm{y-\tilde{y}_{\mathcal{T}}(u)}^2 + \norm{\phi-\tilde{\phi}_{\mathcal{T}}(u)}^2 + \norm{u-\Pi_h^0 u}^2\\& \qquad+(\phi+\lambda u - \Pi_h^0 (\phi+\lambda u), \Pi_h^0 u-u)\lesssim \norm{y-\tilde{y}_{\mathcal{T}}(u)}^2 + \norm{\phi-\tilde{\phi}_{\mathcal{T}}(u)}^2\\&\;\;\; \qquad + \norm{u-\Pi_h^0 u}^2 +\norm{(\phi+\lambda u) - \Pi_h^0 (\phi+\lambda u)} \norm{\Pi_h^0 u-u}.     \end{align*}
The estimates of $\norm{y-\tilde{y}_{\mathcal{T}}(u)}$ and $\norm{\phi-\tilde{\phi}_{\mathcal{T}}(u)}$ follows from Lemmas \ref{lem:estm_aux_state}-\ref{lem:estm_aux_adj_state}. From the regularity of solution from Lemma \ref{lem:regularity_WC}, we obtain the second estimate
$\norm{u-u_h}+\norm{y-y_{\mathcal{T}}}\lesssim h.$ This completes the lemma.   
\end{proof}

\begin{theorem}[Energy error estimate of state]
There holds 
\begin{align*}
    \norm{\widehat{I}_h y - \widehat{y}_h}_{1,h} &\lesssim ( \sum_{T\in \mathcal{T}_h} (\norm{\nabla(y-\mathcal{E}_{\mathcal{T}}^{1}y)}_T^2 + h_T\norm{\nabla(y-\mathcal{E}_{\mathcal{T}}^{1}y)}_{\partial T}^2))^{1/2} + \norm{u-u_h}.
\end{align*}
Furthermore, using the regularity of solution from Lemma \ref{lem:regularity_WC}, we obtain the following convergence rate
$$\norm{\widehat{I}_h y - \widehat{y}_h}_{1,h} = \mathcal{O} (h).$$
\end{theorem}
\begin{proof}
 The proof of the first part follows, similar to the proof of Theorem \ref{thm:energy_estm_state}. For the second part, we use the regularity of solutions from Lemma \ref{lem:regularity_WC}, and we have the estimate $\norm{u-u_h}= \mathcal{O}(h)$ from Lemma \ref{thm:control_conv_rates_WC}, and $\norm{\nabla(y-\mathcal{E}_{\mathcal{T}}^{1}y)}_T= \mathcal{O}(h_T), \norm{\nabla(y-\mathcal{E}_{\mathcal{T}}^{1}y)}_{\partial T}= \mathcal{O}(h_T).$ Hence the second estimate follows.
\end{proof}

\begin{theorem}[Energy error estimate of adjoint state]
There holds 
\begin{align*}
    \norm{\widehat{I}_h \phi - \widehat{\phi}_h}_{1,h} &\lesssim  ( \sum_{T\in \mathcal{T}_h} (\norm{\nabla(\phi-\mathcal{E}_{\mathcal{T}}^{1}\phi)}_T^2 + h_T \norm{\nabla(\phi-\mathcal{E}_{\mathcal{T}}^{1}\phi)}_{\partial T}^2))^{1/2} \\&+ \norm{u-u_h}+\norm{y-\tilde{y}_{\mathcal{T}}(u)}.
\end{align*}
Furthermore, the regularity of solution from Lemma \ref{lem:regularity_WC} leads to
$$\norm{\widehat{I}_h \phi - \widehat{\phi}_h}_{1,h} = \mathcal{O}(h).$$   
\end{theorem}
\begin{proof}
	 The proof of the first part follows, similar to the proof of Theorem \ref{lem:estm_aux_adj_state}. For the second part, we use the regularity of solutions from Lemma \ref{lem:regularity_WC}, and we have the estimate $\norm{u-u_h} = \mathcal{O}(h)$ from Lemma \ref{thm:control_conv_rates_WC}, and estimate $\norm{y-\tilde{y}_{\mathcal{T}}(u)} = \mathcal{O} (h)$ follows from Lemma \ref{lem:estm_aux_state}, $\norm{\nabla(y-\mathcal{E}_{\mathcal{T}}^{1}y)}_T = \mathcal{O} (h_T), \norm{\nabla(y-\mathcal{E}_{\mathcal{T}}^{1}y)}_{\partial T} = \mathcal{O}(h_T).$ Hence the second estimate follows.
\end{proof}

\subsection{Scheme WC2 (Variational discretization):}\label{Scheme:WC_3}
This scheme considers the variational discretization approach. The discrete control problem is as follows:
\begin{equation}\label{eqn:disc_costfun_WCVD}
	\min J(\widehat{y}_h, u)=\frac{1}{2}\lVert y_{\mathcal{T}}-y_{d}\rVert^{2}+\frac{\lambda}{2}\|u\|^{2}, 
\end{equation} 
subject to the PDE,
\begin{equation}\label{eqn:disc_state_WCVD}
	a_h(\widehat{y}_h, \widehat{w}_h)=(f,w_{\mathcal{T}})+(u,w_{\mathcal{T}}) \quad \forall \widehat{w}_h \in \widehat{V}_{h,0}^{1+}.    
\end{equation}
The control $u \in U_{ad}.$ 

\noindent
Existence and uniqueness of the solution of the problem \eqref{eqn:disc_costfun_WCVD}-\eqref{eqn:disc_state_WCVD} follows from Theorem \ref{thm:existance_thm}. The discrete optimality system reads as: find $(\widehat{y}_h, u_h, \widehat{\phi}_h)\in \widehat{V}_{h,0}^{1+} \times U_{ad} \times \widehat{V}_{h,0}^{1+} $ s.t.
\begin{align}
	a_h(\widehat{y}_h, \widehat{w}_h)&=(f, w_{\mathcal{T}})+(u_h, w_{\mathcal{T}}) \quad \forall \widehat{w}_h \in \widehat{V}_{h,0}^{1+}\label{eqn:disc_state_VD},\\
	a_h(\widehat{w}_h, \widehat{\phi}_h)&=(y_{\mathcal{T}}, w_{\mathcal{T}})-(y_d, w_{\mathcal{T}} ) \quad \forall \widehat{w}_h \in \widehat{V}_{h, 0}^{1+}\label{eqn:disc_adjstate_VD},\\
	(\phi_{\mathcal{T}} + \lambda u_h, v-u_h)&\geq 0 \quad \forall v\in U_{ad}.  \label{eqn:disc_VD} 
\end{align}

We need the following auxiliary problems for error analysis: Let $\tilde{y}_h(u), \tilde{\phi}_h(u) \in \widehat{V}_{h, 0}^{1}$ solve the following:
\begin{align}
 a_h(\tilde{y}_h(u), \widehat{w}_h)&=(f, w_{\mathcal{T}} )+(u, w_{\mathcal{T}} ) \quad \forall \widehat{w}_h \in \widehat{V}_{h, 0}^{1},\label{eqn:aux_y(u)_WCM1} \\
a_h(\widehat{w}_h, \tilde{\phi}_h(u))&=(y-y_d, w_{\mathcal{T}} ) \quad \forall \widehat{w}_h \in \widehat{V}_{h, 0}^{1}\label{eqn:aux_p(u)_WCM1}.
\end{align}

\begin{lemma}[$L^2$-error estimate of control and state]\label{thm:control_conv_rates_VD}
Let $(y,u)$ be the optimal state and control of \eqref{eqn:conti_os_state}-\eqref{eqn:conti_os_control}. The pair $(\hhoelem{y},u_h)$ is the discrete state and optimal control \eqref{eqn:disc_state_VD}-\eqref{eqn:disc_VD}.  Then, it holds
	\begin{align*}
		\norm{u-u_h} + \norm{y-y_{\mathcal{T}}} &\lesssim \norm{y - \tilde{y}_{\mathcal{T}}(u)}+ \norm{\phi - \tilde{\phi}_{\mathcal{T}}(u)},
	\end{align*}
where $\tilde{\phi}_h(u), \tilde{y}_h(u)$ solves \eqref{eqn:aux_y(u)_WCM1} and \eqref{eqn:aux_p(u)_WCM1} respectively. Further, the regularity of solution from Lemma \ref{lem:regularity_WC} leads to
\begin{align*}
\norm{u-u_h} + \norm{y-y_{\mathcal{T}}}
= \mathcal{O}(h^3).
\end{align*}
\end{lemma}

\begin{proof}
The choice $ v= u_h$ in \eqref{eqn:conti_VI_WC} leads to  
\begin{align}\label{aux:VD_01}
		(\phi+\lambda u, u_h-u)\geq 0.
	\end{align}
Next, the choice $v=u$ in \eqref{eqn:disc_VD} reveals
\begin{align}\label{aux:VD_02}
(\phi_{\mathcal{T}}+\lambda u_h, u-u_h)\geq 0.
\end{align}
Add \eqref{aux:VD_01} and \eqref{aux:VD_02} to get $(\phi_{\mathcal{T}}-\phi +\lambda(u_h-u), u-u_h)\geq 0.$ The elementary manipulations lead to $(\phi_{\mathcal{T}}-\phi, u-u_h)  - \lambda \norm{u-u_h}^2 \geq 0$. Introduce the split to get $\lambda\norm{u-u_h}^2\leq (\phi_{\mathcal{T}}-\tilde{\phi}_{\mathcal{T}}(u),u-u_h)+(\tilde{\phi}_{\mathcal{T}}(u)-\phi,u-u_h)\leq (\tilde{\phi}_{\mathcal{T}}(u)-\phi,u-u_h)+(y-y_{\mathcal{T}}, \tilde{y}_{\mathcal{T}}(u)-y)-\norm{y-y_{\mathcal{T}}}^2.$ The Cauchy-Schwarz inequality yields 
\begin{align}
\norm{u-u_h} + \norm{y-y_{\mathcal{T}}} &\lesssim \norm{y-\tilde{y}_{\mathcal{T}}(u)}+\norm{\phi -\tilde{\phi}_{\mathcal{T}}(u)}.\nonumber			
\end{align}
This completes the proof of the first part of Lemma. The proof of the second part follows from the estimates $\norm{y-\tilde{y}_{\mathcal{T}}(u)}= \mathcal{O} (h^3)$ and $\norm{\phi-\tilde{\phi}_{\mathcal{T}}(u)}= \mathcal{O} (h^3)$. Those two estimates can be derived using the ideas from the proof of Lemmas \ref{lem:estm_aux_state}-\ref{lem:estm_aux_adj_state}. 
\end{proof}

\noindent
Below, we state the energy estimates of state and adjoint state variables. The proof follows from the proof of Lemmas \ref{thm:energy_estm_state}-\ref{thm:energy_estm_adj_state_S1UC}.
 
\begin{lemma}[Energy error estimate of state]\label{thm:energy_estm_state_VD}
	There holds 
	\begin{align*}
		\norm{\widehat{I}_h y - \widehat{y}_h}_{1,h} &\lesssim ( \sum_{T\in \mathcal{T}_h} (\norm{\nabla(y-\mathcal{E}_{\mathcal{T}}^{2}y)}_T^2 + h_T\norm{\nabla(y-\mathcal{E}_{\mathcal{T}}^{2}y)}_{\partial T}^2))^{1/2} + \norm{u-u_h}.
	\end{align*}
 Further, the regularity of solution from Lemma \ref{lem:regularity_WC} leads to
 $$\norm{\widehat{I}_h y - \widehat{y}_h}_{1,h}=\mathcal{O} (h^2).$$
\end{lemma}
\begin{lemma}[Energy error estimate of adjoint state]\label{thm:energy_estm_adjstate_VD}
	There holds, 
	\begin{align*}
		\norm{\widehat{I}_h \phi - \widehat{\phi}_h}_{1,h} &\lesssim  ( \sum_{T\in \mathcal{T}_h} (\norm{\nabla(\phi-\mathcal{E}_{\mathcal{T}}^{2}\phi)}_T^2 + h_T\norm{\nabla(\phi-\mathcal{E}_{\mathcal{T}}^{2}\phi)}_{\partial T}^2))^{1/2} + \norm{u-u_h}+\norm{y-\tilde{y}_{\mathcal{T}}(u)}.
	\end{align*} 
Further, the regularity of solution from Lemma \ref{lem:regularity_WC} leads to
$$\norm{\widehat{I}_h \phi - \widehat{\phi}_h}_{1,h}=\mathcal{O}(h^2).$$  
\end{lemma}
\section{Numerical Experiments}\label{sec:Numerical_experiments}
\begin{figure}[h!]
    \centering
    \begin{subfigure}[b]{0.45\linewidth}
        \centering
        \includegraphics[width=\linewidth]{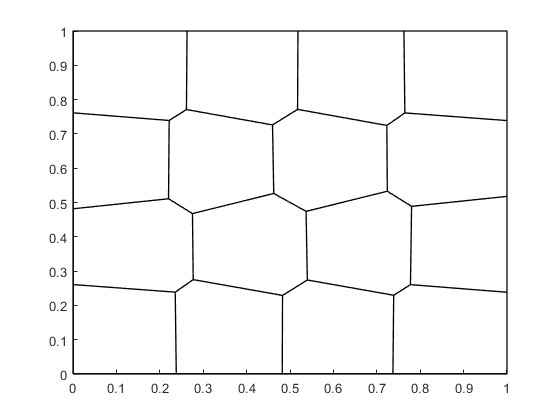}
        \caption{Polygonal Mesh}
        \label{fig:polygonal_mesh}
    \end{subfigure}
    \hfill
    \begin{subfigure}[b]{0.45\linewidth}
        \centering
        \includegraphics[width=\linewidth]{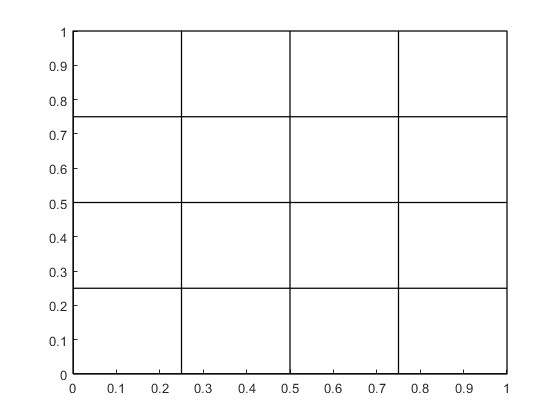}
        \caption{Rectangular Mesh}
        \label{fig:rectangular_mesh}
    \end{subfigure}
    \caption{Meshes}
    \label{fig:combined_meshes}
\end{figure}

In this section, we validate the a priori error estimates
for the error in state, adjoint state, and control variables numerically. For the computations we construct a model problem with known solutions. The convergence rates and robustness of HHO methods are demonstrated. Two $h$-refined mesh families (Rectangular and polygonal (Voronoi-like)) meshes have been used for the computations; see Figure \ref{fig:combined_meshes}.

We consider the domain $\Omega = (0,1)^2.$ The HHO schemes are implemented for the polynomial degree $k\in \{0,1,2,3\}$  and meshes with element numbers $\{16,64,256,1024, 4096\}.$ The errors are measured in $L^2$-norm ($\norm{y-\mathcal{R}\hat{y}_h}$) and energy norm ( $\norm{\mathcal{I}_h y - \hhoelem{y}}_{1,h}$ or  $\norm{\nabla(y-\mathcal{R}\hat{y}_h)})$). 

\begin{table}[h!]
    \centering
    \begin{tabular}{|c|c|c|}
        \hline
        Schemes & Cartesian mesh & polygonal mesh \\ \hline\hline
        Scheme UC1/UC2 &  Figure \ref{fig:scheme_UC1_2_ex2_rectangle} & Figure \ref{fig:scheme_UC1_2_ex2_polygon} \\ \hline
       Scheme UC3.1 & Figure \ref{fig:scheme_UC31_ex1_rectangle} & Figure \ref{fig:scheme_UC31_ex1_polygon} \\ \hline
        Scheme UC3.2 & Figure \ref{fig:scheme_UC32_ex1_rectangle} & Figure \ref{fig:scheme_UC32_ex1_polygon} \\ \hline
         Scheme WC1 & Figure \ref{fig:scheme_WC1_cartesian} & Figure \ref{fig:scheme_WC1_polygonal} \\ \hline
    \end{tabular}
    \caption{Table of references of convergences}
    \label{tab:example}
\end{table}

\subsection{Scheme UC1/UC2:}  The formulations and analyses of Schemes UC1 and UC2 differ in their theoretical aspects; however, their computational implementations yield identical results because there is no control constraint imposed on this problem. In this subsection, we focus on examining the numerical performance of both schemes, UC1 and UC2, to highlight their effectiveness.

To carry out the numerical experiments, we consider the following exact solutions: the state variable \(y(x_1, x_2) = 100\exp(x_1 + x_2)\), the adjoint state \(\phi(x_1, x_2) = \exp(x_1 + x_2) \sin(\pi x_1) \sin(\pi x_2)\), and the control \(u = -\phi / \lambda\). Using these exact solutions, the data functions are specified as \(f = -\Delta y - u\) and \(y_d = y + \Delta y\). For the regularization parameter, we set \(\lambda = 10^{-2}\) to control the balance between the data fidelity and regularization terms.

\begin{figure}[h!]
    \centering
    \setlength{\tabcolsep}{0pt}
    \begin{tabular}{cc}
    \includegraphics[width=0.5\textwidth]{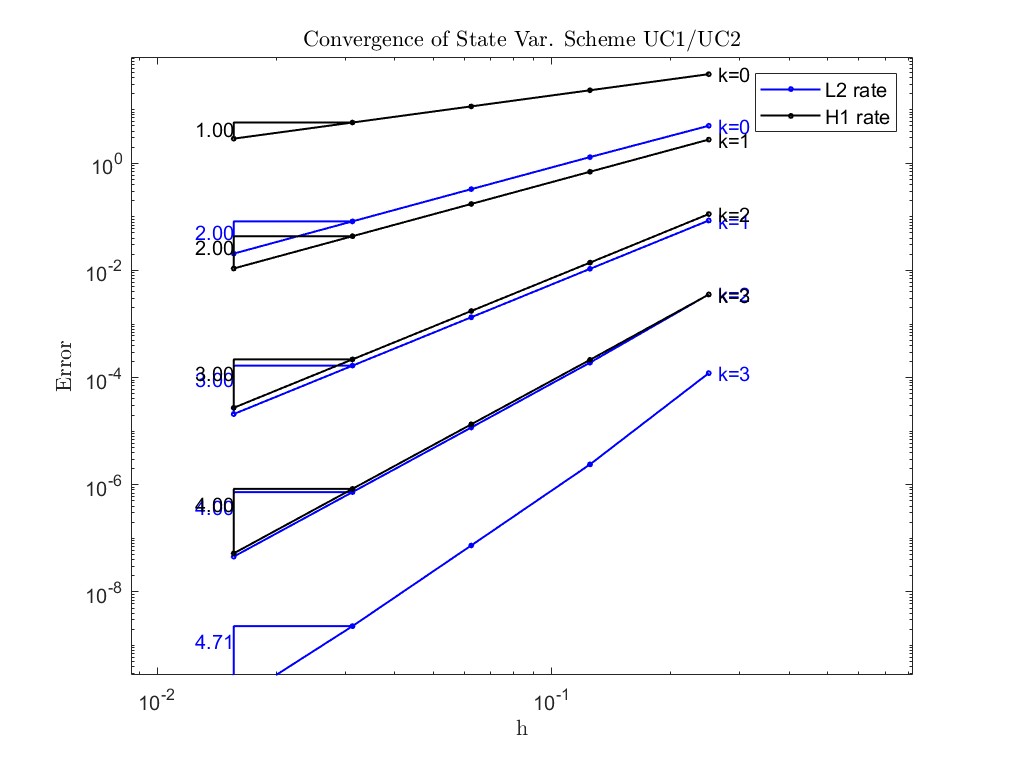}
    & \includegraphics[width=0.5\textwidth]{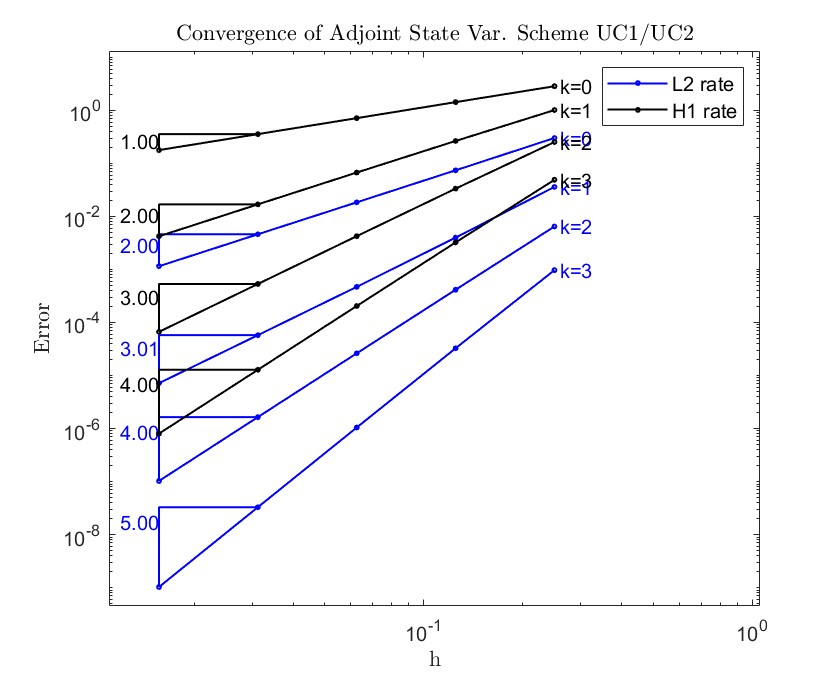}\\ 
    \end{tabular}
\includegraphics[width=0.5\linewidth]{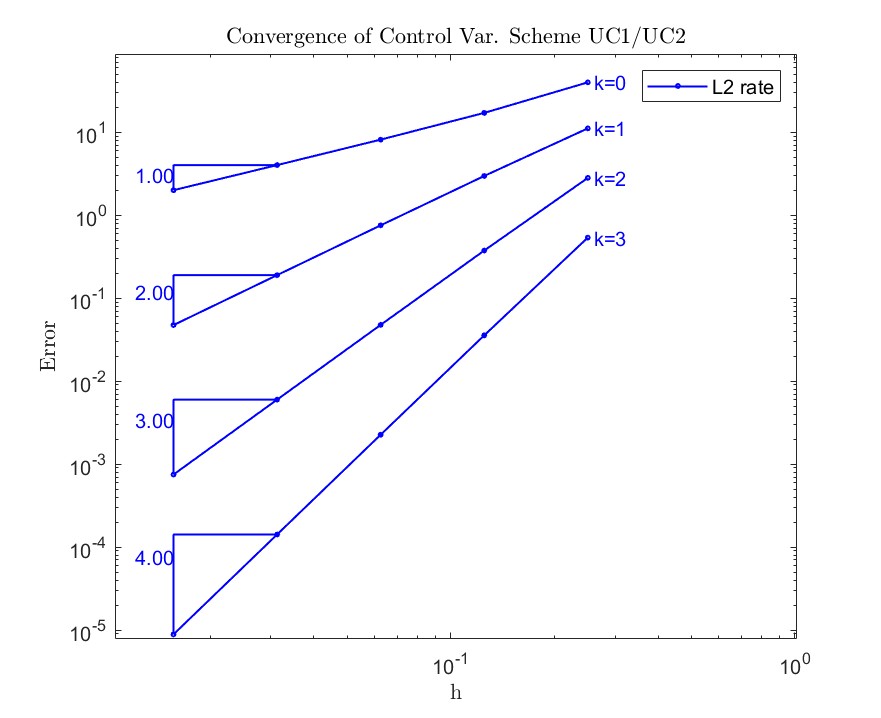}
    \caption{Convergence results for Scheme UC1 \& UC2 on rectangular mesh.}
    \label{fig:scheme_UC1_2_ex2_rectangle}
\end{figure}

\begin{figure}[h!]
    \centering
    \setlength{\tabcolsep}{0pt}
    \begin{tabular}{cc}
    \includegraphics[width=0.5\textwidth]{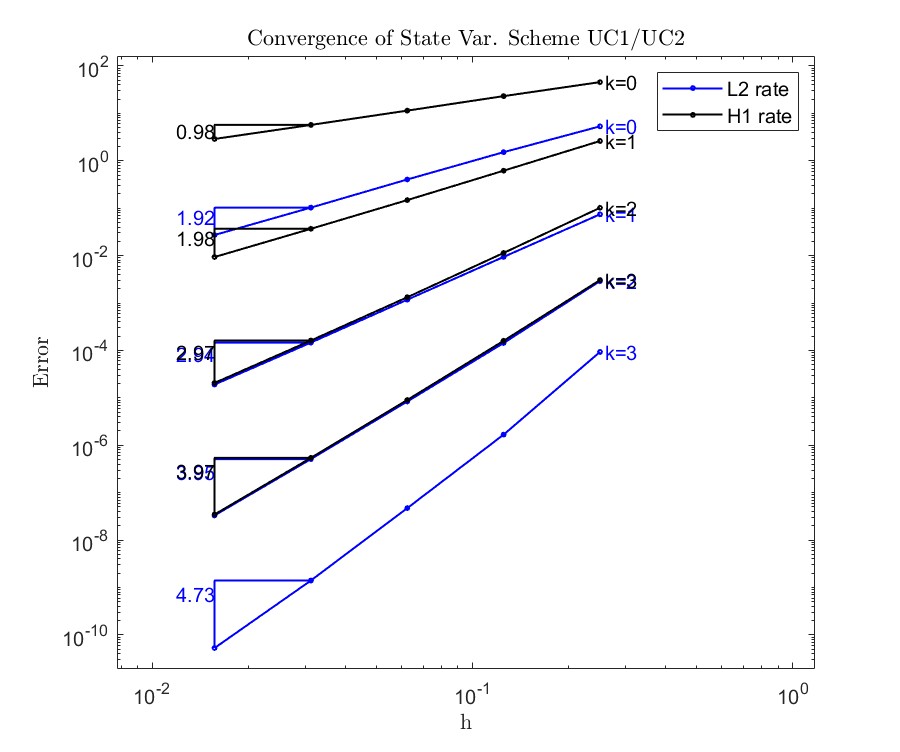}
    & \includegraphics[width=0.5\textwidth]{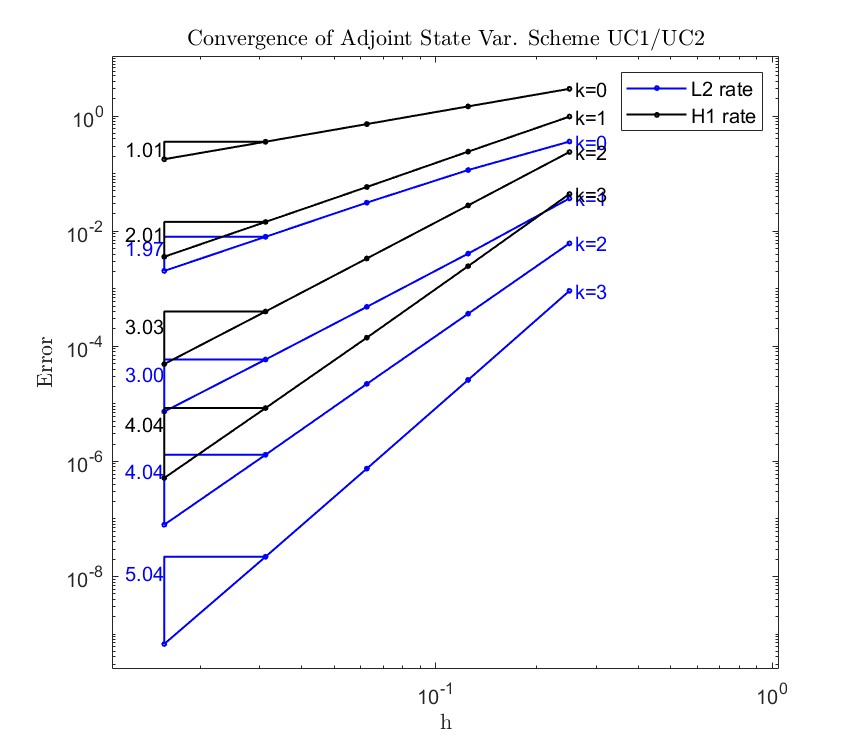}\\ 
    \end{tabular}
\includegraphics[width=0.5\linewidth]{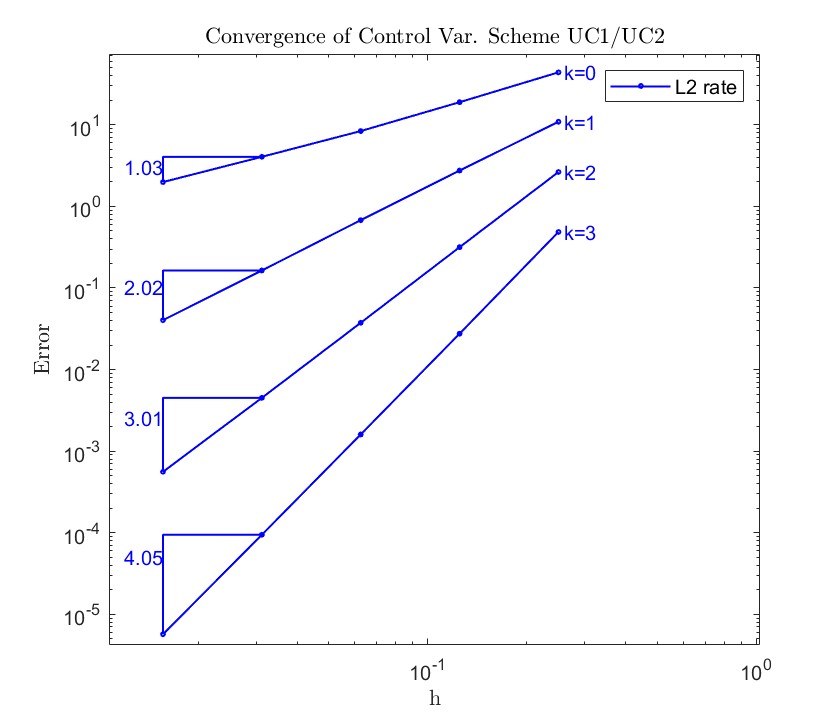}
    \caption{Convergence results for Scheme UC1 \& UC2 on polygonal mesh.}
    \label{fig:scheme_UC1_2_ex2_polygon}
\end{figure}

The convergence of the computed solutions is demonstrated in Figure \ref{fig:scheme_UC1_2_ex2_rectangle} for Cartesian meshes and in Figure \ref{fig:scheme_UC1_2_ex2_polygon} for polygonal meshes. These visualizations highlight the accuracy and reliability of Schemes UC1 and UC2 across different mesh types, illustrating their numerical consistency and robustness.

\subsection{Scheme UC3.1:}
In this subsection, we present the computations for Scheme UC3.1, referred to as the \textit{full reconstruction approach}. In this scheme, the control variable is approximated within a reconstructed HHO space, meaning the control is reconstruction of certain HHO elements. Specifically, the underlying approximation space for the control is reconstruction of zeroth-order HHO space, i.e., piecewise constant approximations over the mesh elements. Therefore the control space consists of piecewise linear polynomials. The state and adjoint state are chosen from the zeroth-order HHO space. 

To perform the numerical experiments, we use the following exact solutions: 
the state variable \(y(x_1, x_2) = 20\sin(2 \pi x_1) \sin(2 \pi x_2)\), 
the adjoint state \(\phi(x_1, x_2) = 5\exp(x_1 + x_2) \sin(\pi x_1) \sin(\pi x_2)\), 
and the control \(u = -\phi / \lambda\). 
Based on these exact solutions, the data functions are defined as \(f = -\Delta y - u\) and \(y_d = y + \Delta y\). The regularization parameter is set to \(\lambda = 10^{-1}\), balancing the trade-off between data fidelity and regularization.

\begin{figure}[h!]
    \centering
    \setlength{\tabcolsep}{0pt}
    \begin{tabular}{cc}
    \includegraphics[width=0.5\textwidth]{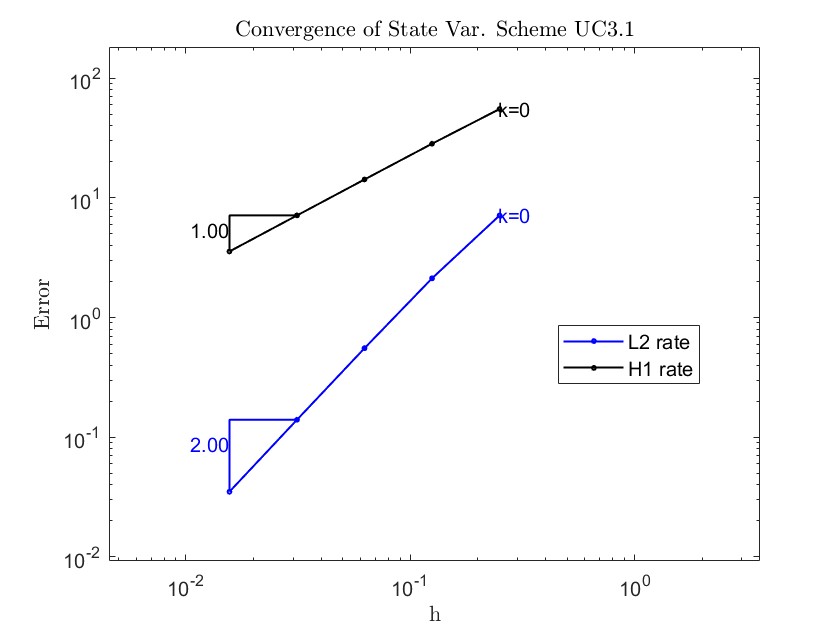}
    & \includegraphics[width=0.5\textwidth]{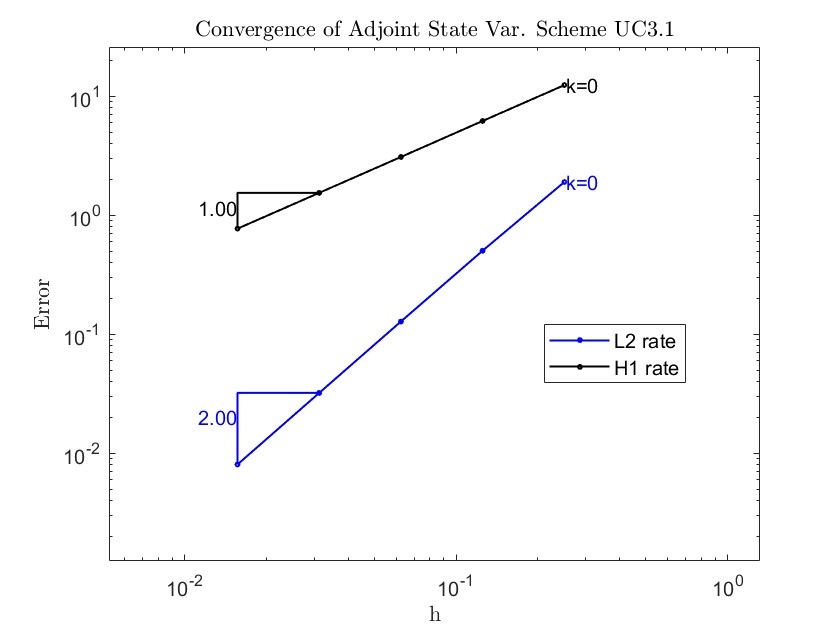}\\ 
    \end{tabular}
\includegraphics[width=0.5\linewidth]{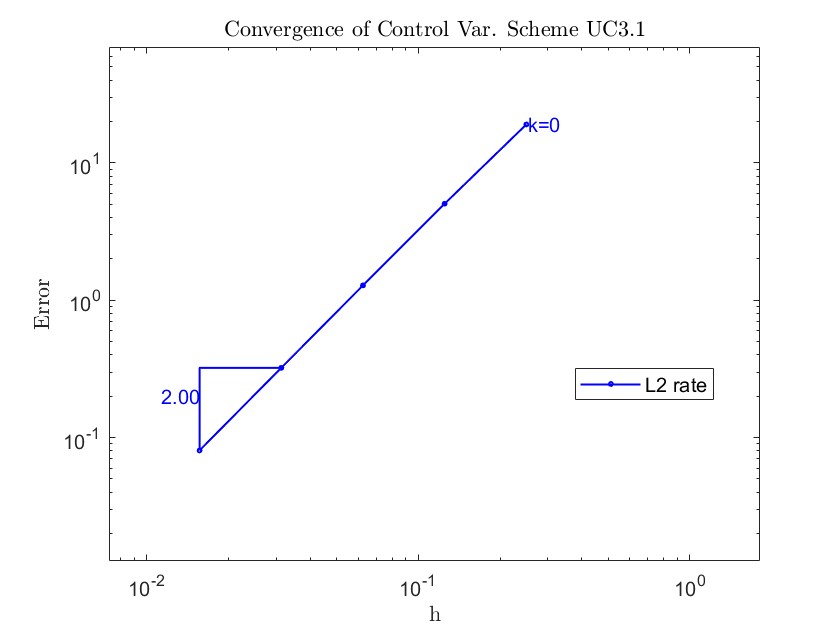}
    \caption{Convergence results for Scheme UC3.1 on rectangular mesh.}
    \label{fig:scheme_UC31_ex1_rectangle}
\end{figure}

\begin{figure}[h!]
    \centering
    \setlength{\tabcolsep}{0pt}
    \begin{tabular}{cc}
    \includegraphics[width=0.5\textwidth]{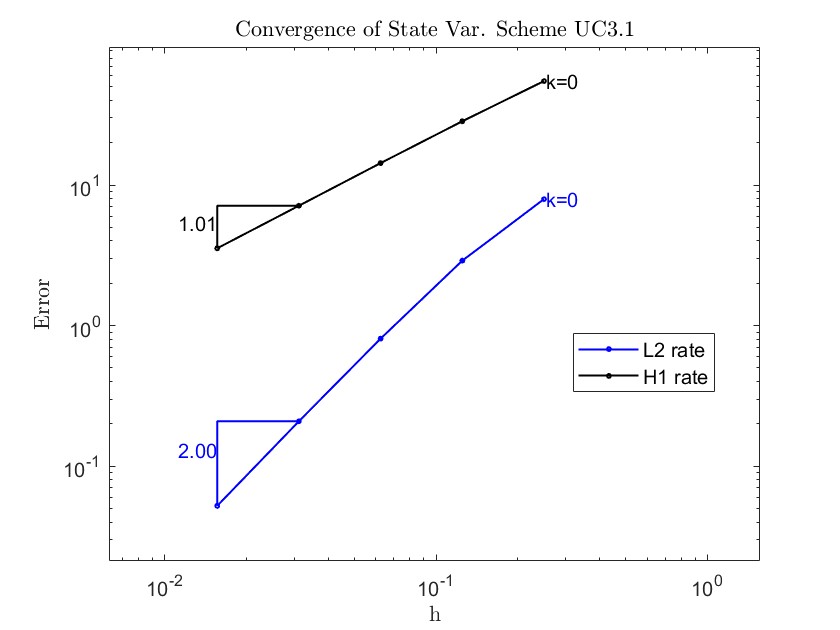}
    & \includegraphics[width=0.5\textwidth]{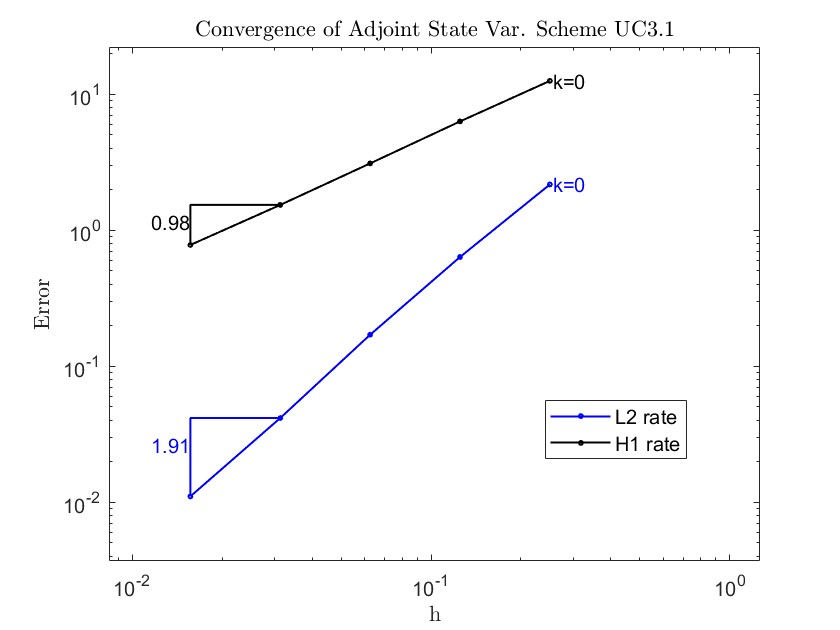}\\ 
    \end{tabular}
\includegraphics[width=0.5\linewidth]{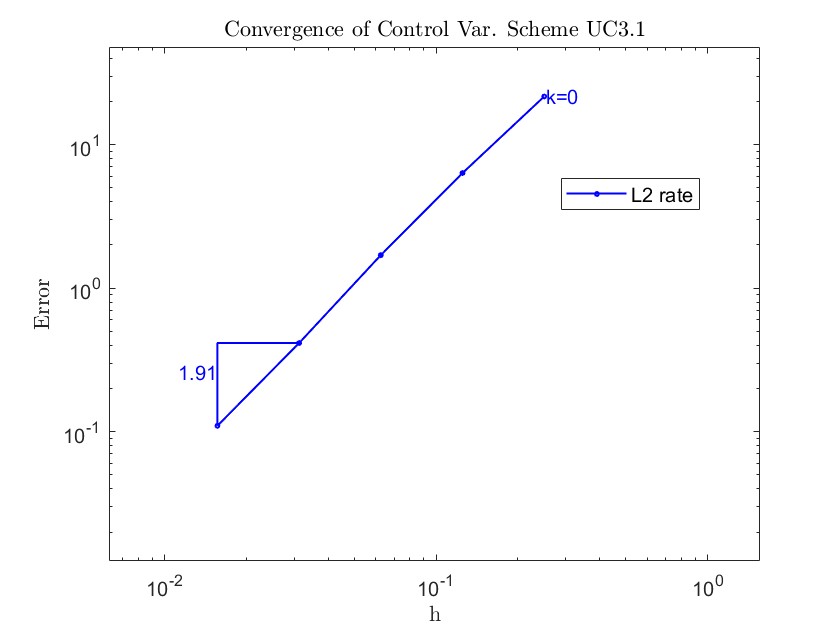}
    \caption{Convergence results for Scheme UC3.1 on polygonal mesh.}
    \label{fig:scheme_UC31_ex1_polygon}
\end{figure}

The convergence of the computed solutions is illustrated in Figure \ref{fig:scheme_UC31_ex1_rectangle} for Cartesian meshes and in Figure \ref{fig:scheme_UC31_ex1_polygon} for polygonal meshes. These figures emphasize the accuracy and robustness of Scheme UC3.1 across different mesh configurations.

\subsection{Scheme UC3.2:}
In this subsection, we present the computations for Scheme UC3.2, referred to as the \textit{partial reconstruction approach}. In this scheme, the control variable is approximated within a reconstructed HHO space, meaning the control is reconstruction of $k$-th order HHO elements. Specifically, the underlying approximation space for the control is reconstruction of $k$-th order HHO space, i.e., piecewise $k$-th order polynomial approximations over the mesh elements. Therefore the control space consists of piecewise $k$-th order polynomials. The state and adjoint state are chosen from the $k$-th order mixed HHO space, i.e., elements are $k+1$-th order polynomial in the interior and $k$-th order polynomial on the skeleton. For the computations of this scheme, we choose $k = 1,2,\; \text{and}\; 3.$

To perform the numerical experiments, we use the following exact solutions: 
the state variable \(y(x_1, x_2) = \sin(2 \pi x_1) \sin(2 \pi x_2)\), 
the adjoint state \(\phi(x_1, x_2) = \exp(x_1 + x_2) \sin(\pi x_1) \sin(\pi x_2)\), 
and the control \(u = -\phi / \lambda\). 
Based on these exact solutions, the data functions are defined as \(f = -\Delta y - u\) and \(y_d = y + \Delta y\). The regularization parameter is set to \(\lambda = 10^{-2}\), balancing the trade-off between data fidelity and regularization.

\begin{figure}[h!]
    \centering
    \setlength{\tabcolsep}{0pt}
    \begin{tabular}{cc}
    \includegraphics[width=0.5\textwidth]{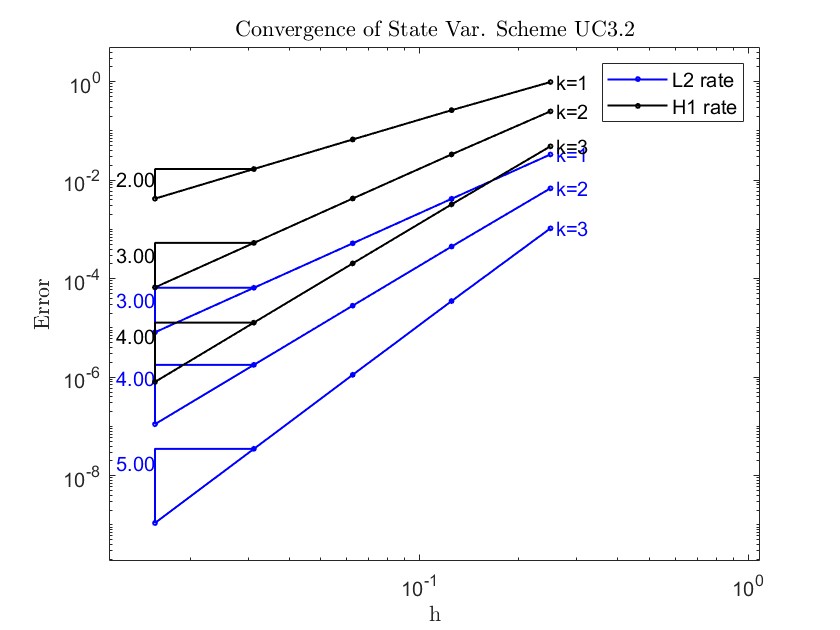}
    & \includegraphics[width=0.5\textwidth]{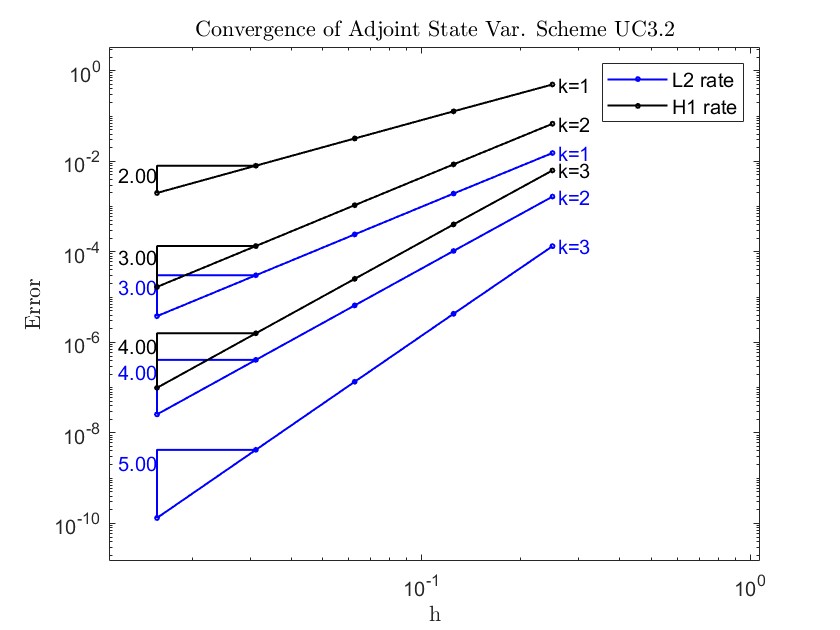}\\ 
    \end{tabular}
\includegraphics[width=0.5\linewidth]{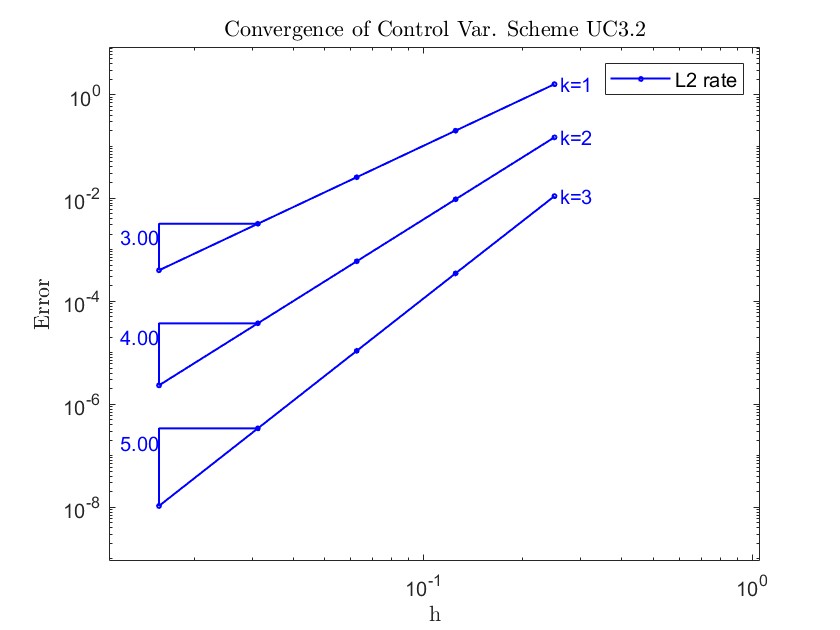}
    \caption{Convergence results for Scheme UC3.2 on rectangular mesh.}
    \label{fig:scheme_UC32_ex1_rectangle}
\end{figure}

\begin{figure}[h!]
    \centering
    \setlength{\tabcolsep}{0pt}
    \begin{tabular}{cc}
    \includegraphics[width=0.5\textwidth]{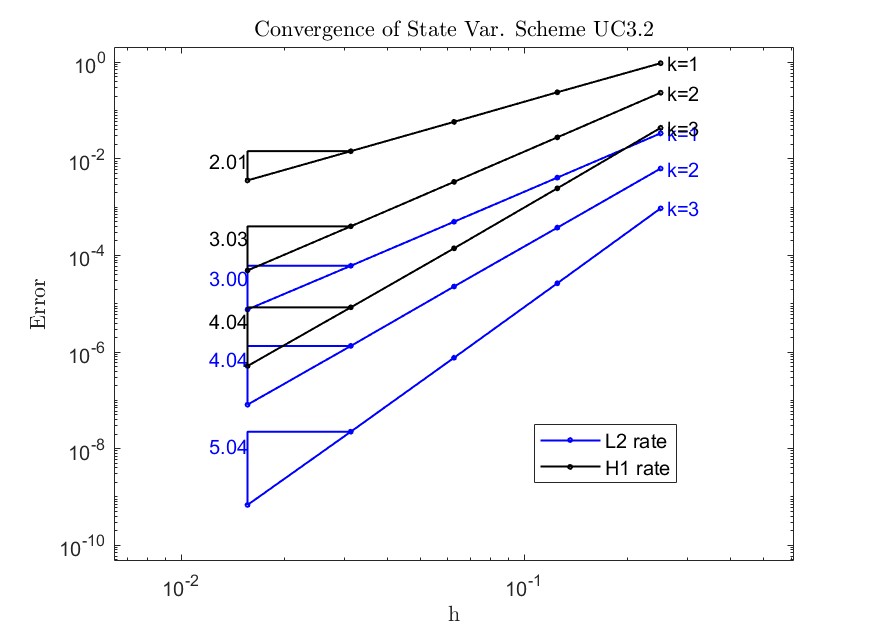}
    & \includegraphics[width=0.5\textwidth]{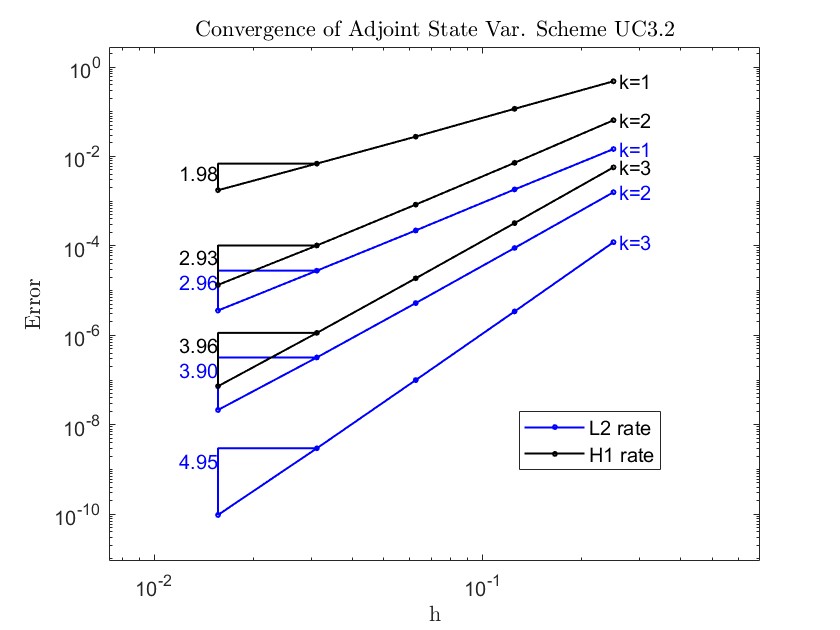}\\ 
    \end{tabular}
\includegraphics[width=0.5\linewidth]{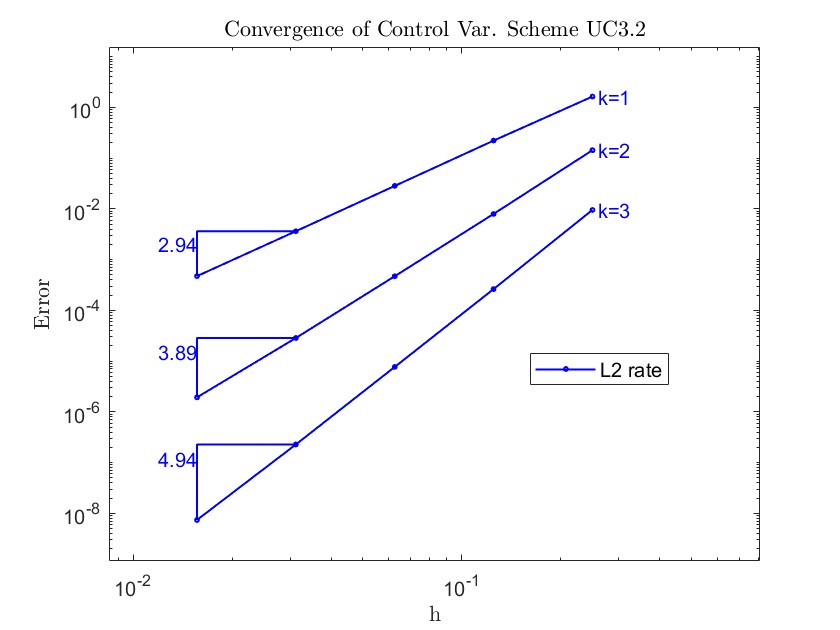}
    \caption{Convergence results for Scheme UC3.2 on polygonal mesh.}
    \label{fig:scheme_UC32_ex1_polygon}
\end{figure}

The convergence of the computed solutions is illustrated in Figure \ref{fig:scheme_UC32_ex1_rectangle} for Cartesian meshes and in Figure \ref{fig:scheme_UC32_ex1_polygon} for polygonal meshes. These figures emphasize the accuracy and robustness of Scheme UC3.2 across different mesh configurations.

\subsection{Scheme WC1:}
In this subsection, we present the computations for Scheme WC1. In this scheme, the control variable is approximated by piecewise constant polynomials. The state and adjoint state variables are approximated by $0$-th order HHO space. The Control satisfies the constraints on each element. Thus it satisfies constraints in the whole domain $\Omega$. We use the projected gradient descent algorithm \cite[Subsection 2.12.2]{trolzstch:2005:Book} to solve the problem \eqref{eqn:disc_state_WCM0}-\eqref{eqn:disc_VI_WCM0}.

To perform the numerical experiments, we use the following exact solutions: 
the state variable \(y(x_1, x_2) = \sin(2 \pi x_1) \sin(2 \pi x_2)\), 
the adjoint state \(\phi(x_1, x_2) = \exp(x_1 + x_2) \sin(\pi x_1) \sin(\pi x_2)\), 
and the control \(u = \mathcal{P}_{U_{ad}}(-\phi / \lambda)\). 
Based on these exact solutions, the data functions are defined as \(f = -\Delta y - u\) and \(y_d = y + \Delta y\). The regularization parameter is set to \(\lambda = 10^{-2}\), balancing the trade-off between data fidelity and regularization. The constants on the box constraints are $U_a=-250$ and $U_b=-10$.

\begin{figure}[h!]
    \centering
    \setlength{\tabcolsep}{0pt}
    \begin{tabular}{cc}
    \includegraphics[width=0.5\textwidth]{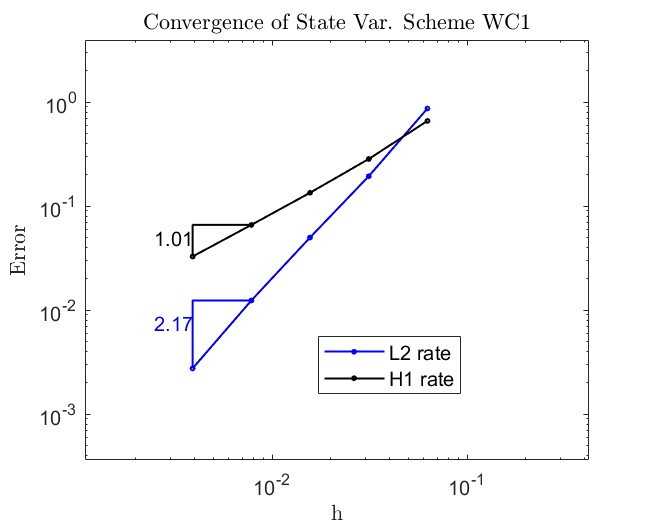}
    & \includegraphics[width=0.5\textwidth]{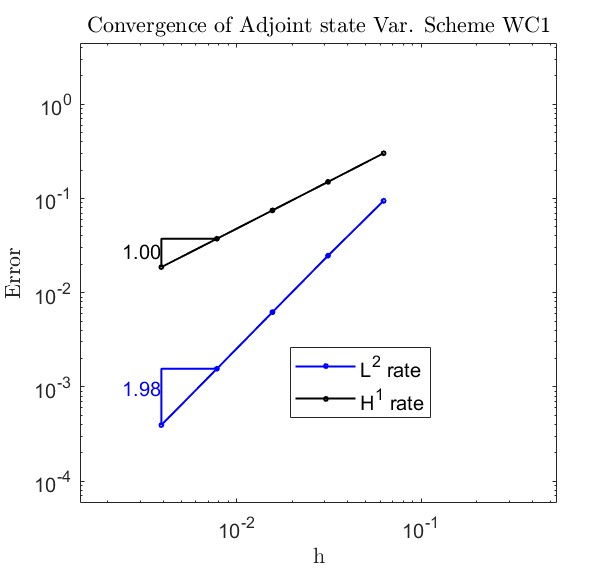}\\ 
    \end{tabular}
\includegraphics[width=0.5\linewidth]{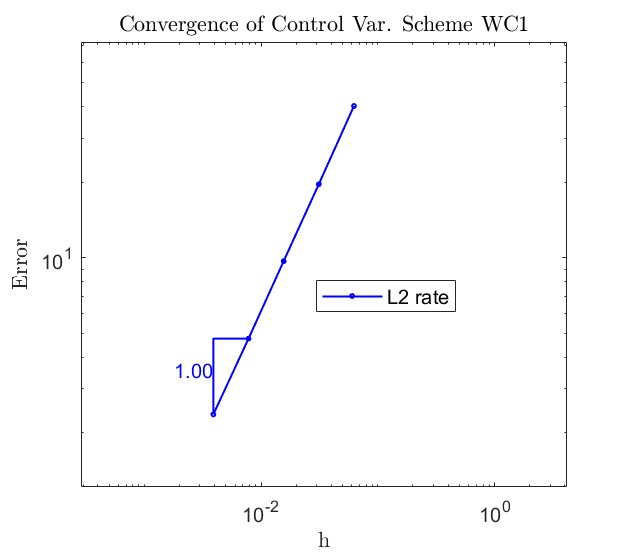}
    \caption{Convergence results for Scheme WC1 on rectangular mesh.}
    \label{fig:scheme_WC1_cartesian}
\end{figure}

\begin{figure}[h!]
    \centering
    \setlength{\tabcolsep}{0pt}
    \begin{tabular}{cc}
    \includegraphics[width=0.5\textwidth]{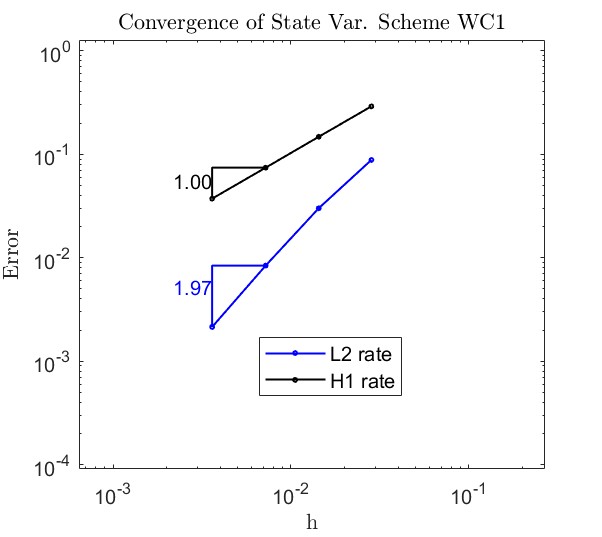}
    & \includegraphics[width=0.5\textwidth]{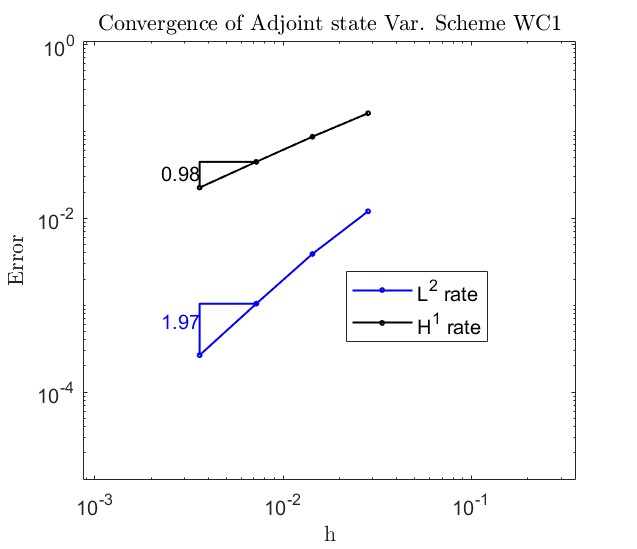}\\ 
    \end{tabular}
\includegraphics[width=0.5\linewidth]{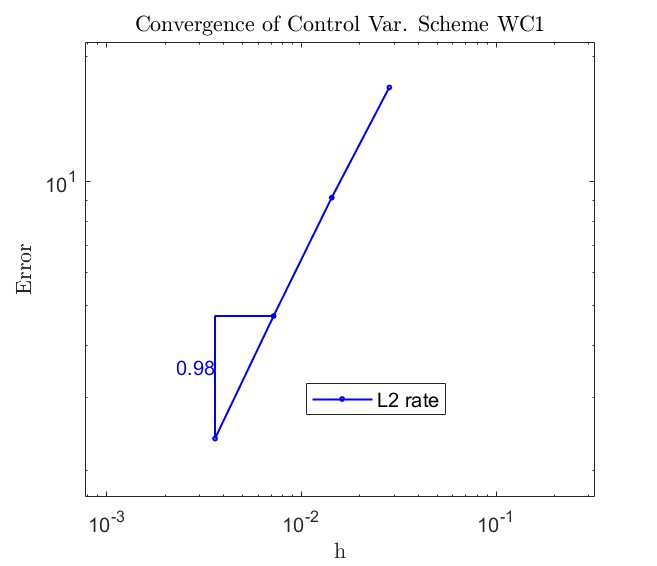}
    \caption{Convergence results for Scheme WC1 on polygonal mesh.}
    \label{fig:scheme_WC1_polygonal}
\end{figure}

The convergence of the computed solutions is illustrated in Figure \ref{fig:scheme_WC1_cartesian} for cartesian meshes and in Figure \ref{fig:scheme_WC1_polygonal} for polygonal meshes. We observe a linear rate of convergence for the $H^1$-error in the state, the adjoint state, and $L^2$- error for the control variables.

\section{Conclusion}
In this paper, we designed and analyzed several Hybrid High-Order approximation schemes to solve distributed optimal control problems governed by the Poisson equation. The study addressed both unconstrained and box-constrained control problems.

For the unconstrained control problem, three novel schemes were proposed. Our methods demonstrate a significant improvement over classical finite element methods. Specifically:
\begin{itemize}
    \item \textit{Enhanced Convergence Rates:} While traditional FEM typically achieves a convergence rate of \( k+1 \) (with \( k \) being the polynomial degree of approximation), our third scheme (UC 3.1\&3.2) i.e., the partial and full reconstruction approaches achieves a convergence rate of \( k+2 \). This improvement is realized by selecting the control variable from a carefully constructed reconstruction space.
    \item \textit{Innovative Approaches:} The methods include full and partial reconstruction strategies that leverage the flexibility of HHO spaces, ensuring better approximation of control variables.
\end{itemize}

For the box-constrained control problem, we introduced two schemes that leverage the unique features of HHO spaces:
\begin{itemize}
    \item \textit{Efficient Use of Low-Order Elements:} The first scheme (WC 1) achieves linear convergence using only the lowest-order elements (\( \mathbb{P}_0 \)) for all variables. This contrasts with FEM, which requires higher-order linear elements for similar results.
    \item \textit{Cubic Convergence:} The second scheme (WC 2) achieves a remarkable cubic convergence for the control variable, provided sufficient regularity assumptions on the data (\( f, u_d \in H^1(\Omega) \)) are met.
\end{itemize}

The theoretical findings were further substantiated through comprehensive numerical experiments, confirming the accuracy and efficiency of our proposed methods. These results underline the potential of HHO methods to serve as powerful tools in solving optimal control problems, offering higher accuracy and flexibility compared to conventional approaches.

In summary, this work presents a robust framework that not only advances the state-of-the-art numerical methods for optimal control but also opens new avenues for further research in applying high-order schemes to complex control problems such as semi-linear, quasi-linear, and fully nonlinear control problems.

\bigskip

\noindent\textbf{Acknowledgements}\\
The first author acknowledges the DST-SERB-MATRICS grant MTR/2023/000681 for financial support and the DST FIST grant SR/FST/MS II/2023/139 for partial financial support. The second author extends gratitude to Prof. Neela Nataraj (Indian Institute of Technology, Bombay, India) and Dr. Raman Kumar (University of Montpellier, France) for their valuable and insightful discussions.

\medskip

\bibliographystyle{plain}
\bibliography{Control}
 \end{document}